\documentclass[12pt]{article}

\usepackage{amsmath,amssymb}
\usepackage{mathrsfs}
\usepackage{natbib}
\RequirePackage[colorlinks,citecolor=blue,urlcolor=blue]{hyperref}
\usepackage{graphicx}
\usepackage{multirow}
\usepackage{float}
\usepackage{url} 
\usepackage{enumerate}
\usepackage{color}
\usepackage{booktabs}
\usepackage{times}
\usepackage{subfig} 
\usepackage{overpic} 
\usepackage{dsfont}
\usepackage{enumitem}
\allowdisplaybreaks[4]

\def\T{{ \mathrm{\scriptscriptstyle T} }}
\DeclareMathOperator*{\argmin}{argmin}

\DeclareMathOperator*{\diam}{diam}
\newcommand*{\mc}[1]{\mathcal{#1}}

\newcommand*{\ms}[1]{\mathsf{#1}}

\newcommand{\E}{\mathbb{E}}
\newcommand{\ind}{\mathds{1}}
\newcommand{\dl}{\mathrm{d}}
\newcommand{\lb}{\left(} 
\newcommand{\rb}{\right)} 
\newcommand{\br}[1]{\lb #1 \rb} 
\newcommand{\abs}[1]{\left| #1 \right|} 

\newtheorem{theorem}{Theorem}[section]

\newtheorem{remark}[theorem]{Remark}
\newtheorem{definition}[theorem]{Definition}

\newcommand{\blind}{0}

\begin{document}

	\def\spacingset#1{\renewcommand{\baselinestretch}%
		{#1}\small\normalsize} \spacingset{1}

	
	\if0\blind
	{
		\title{\bf Semi-supervised Fr\'echet Regression}
		\author{Rui Qiu, Zhou Yu\\
			School of Statistics, East China Normal University\\
			and \\
			Zhenhua Lin \\
			Department of Statistics and Data Science, National University of Singapore}
		\maketitle
	} \fi
	
	\if1\blind
	{
		\bigskip
		\bigskip
		\bigskip
		\begin{center}
			{\LARGE\bf Semi-supervised Fr\'echet Regression}
		\end{center}
		\medskip
	} \fi
	
	\bigskip
	\begin{abstract}
		This paper explores the field of semi-supervised Fr\'echet regression, driven by the significant costs associated with obtaining non-Euclidean labels. Methodologically, we propose two novel methods: semi-supervised NW Fr\'echet regression and semi-supervised kNN Fr\'echet regression, both based on graph distance acquired from all feature instances. These methods extend the scope of existing semi-supervised Euclidean regression methods.  We establish their convergence rates with limited labeled data and large amounts of unlabeled data, taking into account the low-dimensional manifold structure of the feature space. Through comprehensive simulations across diverse settings and applications to real data, we demonstrate the superior performance of our methods over their supervised counterparts. This study addresses existing research gaps and paves the way for further exploration and advancements in the field of semi-supervised Fr\'echet regression.
	\end{abstract}
	
	\noindent%
	\vfill
	
	\newpage
	\spacingset{1.9} 
	\section{Introduction}\label{Introduction}
	Motivated by the demands of contemporary applications, in recent years there has been a growing focus on the statistical analysis of data residing in nonlinear metric spaces. Examples include distributional data in Wasserstein space \citep{delicado2017choosing, chen2023wasserstein}, symmetric positive definite matrices \citep{zhu2009intrinsic, yuan2012local, lin2023additive}, phylogenetic trees \citep{billera2001geometry}, spherical data \citep{di2014nonparametric}, and Riemannian manifold objects \citep{lin2017extrinsic, cornea2017regression}, among others. In regression analysis where a Euclidean predictor is paired with a response situated in a general metric space, which is the focus of this article,  \cite{petersen2019frechet} and \cite{chen2022uniform} developed the methodology and asymptotic theory for both global and local Fr\'echet regression analysis. \cite{schotz2022nonparametric} further studied nonasymptotic behaviors of local Fr\'echet regression, contributing valuable insights into its performance with finite samples. Acquiring non-Euclidean labels, as opposed to Euclidean labels, often demands more labor and resources. Hence, there is a pertinent interest in developing a suitable semi-supervised method for the Fr\'echet regression problem when there is a limited amount of labeled data but a substantial volume of unlabeled data. However, few works consider semi-supervised Fr\'echet regression. 
	
	For data with an Euclidean response, classical semi-supervised learning has become a pivotal research area in machine learning, aiming to address limitations arising from insufficient labeled data in regression or classification tasks. Supervised algorithms necessitate abundant labeled instances to accurately capture complex relationships between variables, but acquiring such data can be arduous, time-consuming, or even impractical in certain domains. Semi-supervised techniques are proposed to overcome these challenges by harnessing the power of both labeled and unlabeled data, thereby improving prediction accuracy. There are significant advancements in this field in classification  \citep[e.g.][]{blum1998combining, ando2005framework, wang2007large, rigollet2007generalization, wang2009efficient, sinha2009semi,wasserman2007statistical, singh2008unlabeled, niyogi2013manifold, gopfert2019can}, parametric regression \citep[e.g.][]{chakrabortty2018efficient, azriel2022semi, rajaraman2022semi, livne2022improved,song2023general}, and nonparametric regression \citep[e.g.][]{belkin2006manifold, wasserman2007statistical, martin2013density}. However, none of these consider response variables sampled from a metric space.
	
	Recent advancements in Euclidean nonparametric semi-supervised regression often hinge upon two fundamental assumptions. The first, known as the manifold assumption, posits that the distribution of the features, denoted as $X$, resides on a low-dimensional manifold $\mathcal{M}$. The second assumption termed the smoothness assumption asserts that the neighboring instances on the manifold tend to possess similar response values $Y$. Building upon these foundational assumptions, semi-supervised regression may be achieved by introducing a regularization term into the objective function based on the samples of $X$  \citep[e.g.,][]{zhu2003semi,belkin2006manifold,niyogi2013manifold}. Alternatively, semi-supervised regression can be achieved by employing unsupervised manifold information extraction techniques, more specifically, via creating a graphical representation of $X$ samples to capture their intrinsic geometrical characteristics. Along this line of research,  \cite{moscovich2017minimax} proposed a geodesic kNN regression method, where the geodesic distances are estimated through shortest-path distances in a graph constructed from all $X$ samples. \cite{belkin2004semi} and \cite{Ji2012ASA}, among others, proposed to encode graph nodes as lower-dimensional vectors while preserving crucial information about their positions and local neighborhood structures. Such graph embedding may be obtained based on Laplace--Beltrami operator on differentiable functions \citep{belkin2004semi}, integral operator on reproducing kernel Hilbert space \citep{Ji2012ASA}, and other methods such as the local linear embedding \citep[e.g.,][]{roweis2000nonlinear} and neural networks \citep{cao2016deep, scarselli2008graph}. 
	
	In this paper, we present kNN type and kernel type of semi-supervised Fr\'echet regression methods and investigate their nonasymptotic risk. Due to the clear geometric interpretation, favorable portability, and convenient calculation means of graph distance, our development focuses on a route based on graph distance. A simple description is that the geodesic distance on a low-dimensional manifold can be efficiently approximated by the length of the shortest-path on a dense feature point graph, according to \cite{tenenbaum2000global}. A large amount of unlabeled data provides us with the possibility of constructing dense graphs, which ensures that we can use graph distances instead of simple Euclidean distances to achieve effective local Fr\'echet regression on a low-dimensional manifold.
	
	There are several advantages in our choice of semi-supervised Fr\'echet methods.
	First, NW regression and kNN regression are two classical local averaging methods, which are simple, effective, interpretable, and widely used in a variety of scenarios, particularly when the response is non-Euclidean. Moreover, many efficient algorithms are available for calculating graph distances, providing great convenience for the implementation of our semi-supervised process. Second, even if the low-dimensional manifold assumption does not hold in reality, the performance of our graph distance based method is comparable to that of NW estimation and kNN regression based on Euclidean distance. Third, we can obtain convergence rate results similar to those under classical Euclidean regression. Notably, our nonasymptotic analysis reveals that, under the low-dimensional manifold assumption as well as limited labeled samples, a substantial amount of unlabeled data can help our method adapt to the intrinsic dimension of the manifold. This offers theoretical guarantees regarding the practical performance of our methods. However, a limitation of our graph distance based semi-supervised method is that generalizing NW Fr\'echet regression to local linear Fr\'echet regression \citep{petersen2019frechet} is difficult. This is mainly because the linear form of the first-order term in the Taylor expansion in local linear regression is no longer appropriate under the manifold assumption. Investigating this aspect further could be a promising direction for future research.
	
	The rest of the paper is organized as follows. In Section~\ref{Semi-Supervised}, we provide an overview of Fr\'echet regression and introduce two noteworthy semi-supervised methods: semi-supervised NW Fr\'echet regression and semi-supervised kNN Fr\'echet regression. Then the convergence rates are established for both methods in Section~\ref{Convergence rate}. Section~\ref{Simulation} encompasses comprehensive simulation studies to compare the performance of two methods under both supervised and semi-supervised learning scenarios. In Section~\ref{Real data}, we showcase the practical merits of semi-supervised Fr\'echet regression through face data. Section~\ref{Discussion} concludes the paper with some discussions. All proofs and additional materials are presented in the supplementary materials.  
	
	\textbf{Notations:} Throughout the paper, we use $c$ or $c'$ to represent general absolute constants whose value may change from line to line. If the value depends on some variables, we indicate them by an index. $S(x, d, \delta)$ denotes the ball of centre $x$  and radius $\delta$ with respect to the metric $d$.  The notation $\hat{m}_x$ uniformly represents the $\hat{m}^{\textsl{snw}}_x, \hat{m}^{\textsl{skn}}, \hat{m}^{\textsl{nw}}$ and $\hat{m}^{\textsl{kn}}$, the four estimation methods which will be introduced successively in the subsequent sections. The specific method referred to by $\hat{m}_x$ can be discerned from the surrounding context.
	
	\section{Semi-supervised Fr\'echet regression}\label{Semi-Supervised}
	Let $(\Omega, d_{\Omega})$ denote a metric space with a specific metric $d_{\Omega}$, and let $\mathcal{R}^{p}$ represent the $p$-dimensional Euclidean space. We consider a random pair $(X, Y) \sim \nu$, where $X \in [0,1]^p \subset \mathcal{R}^{p}$, $Y \in \Omega$, follwing the joint distribution $\nu$.  The marginal distributions of $X$ and $Y$ are denoted as $\nu_X$ and $\nu_Y$, respectively. The conditional distributions $\nu_{X \mid Y}$ and $\nu_{Y \mid X}$ are also assumed to exist.  Our target is to model the regression relationship between $Y$ and $X$. In the case of a general metric space, the concepts of Fr\'echet mean and Fr\'echet variance \citep{frechet1948elements} can be defined as follows
	$$
	\mu_Y=\underset{\omega \in \Omega}{\argmin}\, \E \left\{d_{\Omega}^{2}(Y,\omega)\right\}, \quad V_Y=\E \left\{d_{\Omega}^{2}(Y,\mu_{Y})\right\}.
	$$
	When $\Omega=\mathcal{R}$ and the Euclidean distance is used, $\mu_Y$ and $V_Y$ coincide with the classical mean and variance. The central focus of Fr\'echet regression is to estimate the conditional Fr\'echet mean  \citep{petersen2019frechet}
	\begin{align}\label{goal}
		m_x=\underset{\omega \in \Omega}{\argmin}\, M (x,\omega)=\underset{\omega \in \Omega}{\argmin} \, \E\left\{d_{\Omega}^{2}(Y, \omega) \mid X=x\right\}.
	\end{align}
	Similarly, replacing the intrinsic metric $d_{\Omega}$ of $\Omega$ with Euclidean distance, we return to the classical regression equation
	$$
	m_x=\underset{\omega \in \mathcal{R}}{\argmin}\, \E\left\{(Y-\omega) ^{2}\mid X=x\right\}=\E \big(Y \mid X=x\big).
	$$
	
	The success of semi-supervised Fr\'echet regression, akin to its Euclidean counterpart, hinges upon two critical factors. Firstly, it relies on the concentration of $X$ in proximity to or within lower-dimensional sets (the manifold assumption). Secondly, conditional Fr\'echet mean $m_x$ ought to be related to the marginal distribution of $X$ (the smoothness assumption). Now suppose we have labeled data 
	$$\mathcal{L}_n=\left\{\left(X_1, Y_1\right), \ldots,\left(X_n, Y_n\right)\right\}$$
	from the distribution $\nu$, where the marginal distribution $\nu_X$ is supported on an unknown low-dimensional smooth submanifold $\mathcal{M}$ embedded in $\mathcal{R}^p$,  and $Y_i \in \Omega, 1 \leq i \leq n$. Further, we possess an additional set of unlabeled data
	$$\mathcal{U}_m=\left\{X_{n+1}, \ldots, X_N\right\}$$ 
	from the same distribution but lacking the corresponding $Y$ values. Let $m=N-n$ be the size of $\mathcal{U}_m$.  Compared to non-Euclidean labels, the sample of Euclidean features is easily accessible. Therefore, we assume $m \gg n$.
	
	Noting the manifold structure in $\mathcal{M}$, it is sometimes no longer appropriate to employ the Euclidean distance as a metric for $X$ on the ambient space $\mathcal{R}^p$. Instead, the intrinsic geodesic distance (denoted by $d_{\mathcal{M}}$), which encapsulates valuable information about the underlying manifold shape, serves as a more suitable choice for the feature space. Further, the smoothness assumption enables the utilization of unlabeled data, leading to superior predictions for responses compared to relying solely on labeled data. However, the ensuing difficulty is that the geodesic distance is unknown and it needs to be estimated from available samples. A notable solution called Isomap, proposed by \cite{tenenbaum2000global}, estimates the geodesic distance on $\mathcal{M}$ between data points by the graph distance derived from a constructed graph. The procedure can be summarized as follows:
	
	\begin{itemize}
		\item[1.] Construct an undirected graph $G$ whose vertices represent $\{X_i\}_{i=1}^N$ from labeled and unlabeled data. Pair of points $X_i, X_j, \, 1\leq i \neq j\leq N $ is then connected if and only if $\|X_i - X_j\|_2 \leq r$ based on the Euclidean metric, where $r$ is a predetermined constant. We call the graph constructed in this way $r$-graph.
		\item[2.] For a given instance $x$, set it as a new vertice in the above constructed graph $G$,  and connect it with all vertices in the ball $S(x, \|\cdot\|_2, r)$ of radius $r$.  Compute the shortest-path graph distance $d_G\left(x, X_i\right)$ for all $X_i$ from labeled data. Specifically, the graph distance between $x, X_i$ is defined as
		\begin{align}\label{graph distance}
			d_G(x, X_i)=\min _{\gamma}\left(\left\|x_0-x_1\right\|_2+\ldots+\left\|x_{l-1}-x_l\right\|_2\right),
		\end{align}
		where $P=\left(x_0, \ldots, x_l\right)$ varies over all paths along the edges of $G$ connecting $x\left(=x_0\right)$ to $X_i \left(=x_l\right)$. Actually, the graph distance between any two points from $\{X_i\}_{i=1}^N \cup x$ can be defined similarly.
	\end{itemize}
	
	\begin{figure}[ht!]
		\centering
		\includegraphics[width=0.6\linewidth]{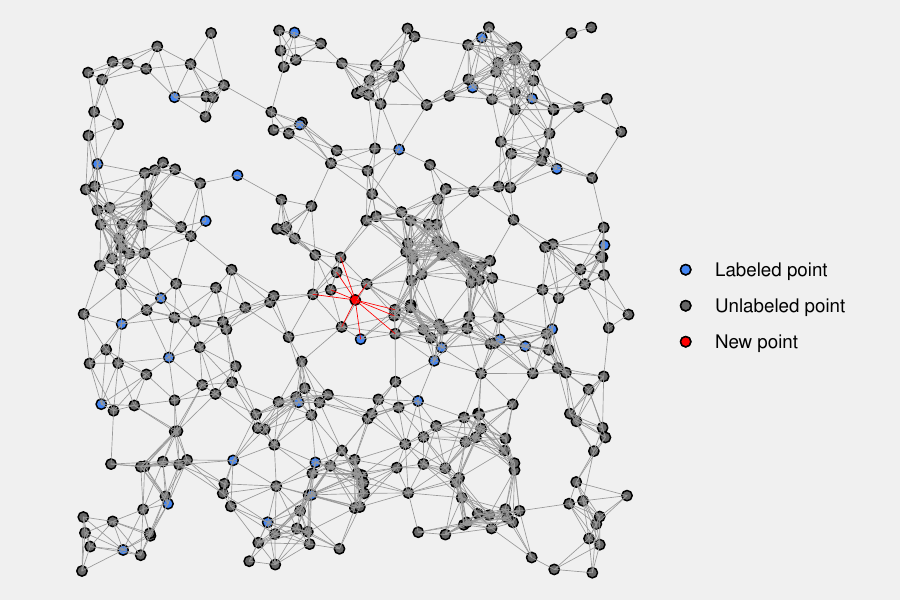}
		\caption{Sparse $r$-graph constructed by $300$ unlabeled points, $30$ labeled points and one new point from $[0, 1]^2$ with $r=0.1$.}
		\label{fig:graph construction}
	\end{figure}
	
	Different from transductive learning, which solely concerns with obtaining label predictions for the given unlabeled data points, inductive learning builds a regressor capable of generating predictions for any object within the feature space. It is worth noting that graph-based methods typically provide an approximation of the geodesic distance only for vertices in the graph. Therefore, in the second step above, any new point is added to the graph in order to implement inductive learning. The whole process is described in Figure~\ref{fig:graph construction}.
	
	\begin{remark}
		Here we consider the utilization of the shortest Euclidean path of $r$-graph to approximate the geodesic distance on the manifold. Additionally, we can consider alternative density-based metrics like \cite{martin2013density}. For instance, we can generalize the geodesic distance to the Fermat distance \citep{hwang2016shortest, fernandez2023intrinsic}, defined as 
		$$
		d_{\mathcal{M}, s}(x, y)=\inf _\gamma \int_I \frac{1}{f\left(\gamma_t\right)^{(s-1) / q}}\|\dot{\gamma}_t\|_2 \dl t,
		$$
		where $x, y \in \mathcal{M}, s \geq 1$, $f(\cdot)$ is the density function on the $q$-dimensional Riemannian manifold $\mathcal{M}$ embedded in $\mathcal{R}^{p}$ and the infimum is taken over all piecewise smooth curves $\gamma: I=[0,1] \rightarrow \mathcal{M}$ with $\gamma(0)=x$, and $\gamma(1)=y$.
		In the special case where $f$ is uniform or $s=1$, the Fermat distance simplifies to (a multiple of) the geodesic distance $d_{\mathcal{M}}$.  Given a sample $\{X_i\}_{i=1}^{n}$, the Fermat distance between $x$ and $y$ can be estimated by
		$$
		d_{n, k}(x, y)=\inf _\gamma \sum_{i=0}^{l-1}\|x_{i+1}-x_i\|_2^s
		$$
		where the infimum is taken over all paths $\gamma=\left(x_0, x_1, \ldots, x_{l}\right)$ with $x_0=x, x_{l}=y$ and $\left\{x_1, x_2, \ldots, x_{l-1}\right\} \subseteq \{X_i\}_{i=1}^{n}$.
		Nevertheless, we argue that the continuity of the regression function with respect to the geodesic distance is more intuitive than with respect to the density-based metric. Hence, we primarily focus on the geodesic distance $d_{\mathcal{M}}$ in our subsequent investigations.
	\end{remark}
	
	\begin{remark}
		Isomap \citep{tenenbaum2000global} stands as a prominent nonlinear dimensionality reduction technique, which also commences by calculating the shortest distances within the constructed $r$-graph mentioned above. To augment this foundation, the authors introduce an additional criterion termed kNN-construction. Every point is connected by an edge to its $k$ nearest neighbors and vice versa.  The resultant graph is denoted as a kNN-graph. We can similarly apply this graph to the subsequent semi-supervised Fr\'echet regression task, as evidenced by the utilization of 4NN-graphs in the real data analysis below. But for achieving inductive learning, unlike the $r$-graph paradigm, the introduction of new sample points might alter the existing connections of edges in the previously constructed kNN-graph due to the kNN-connection criterion. This leads to the fact that for each new sample point, the corresponding kNN-graph needs to be completely updated rather than simply expanded. Leaving aside the theory, we suggest that in practice, a new sample point merely expands the original kNN-graph by connecting it to its $k$ nearest neighbors. And subsequent steps align with those of semi-supervised Fr\'echet regression based on $r$-graphs.
	\end{remark}
	
	\cite{arias2019unconstrained} demonstrates that the graph distance $d_{G}$ and the geodesic distance $d_{\mathcal{M}}$ exhibit remarkable proximity to each other if the data size is large enough, which can be guaranteed by letting $m \to \infty$. This observation allows us to employ the smoothness assumption of $m_x$ with respect to $d_{\mathcal{M}}$ in constructing the local nonparametric Fr\'echet regression based on $d_G$.  For any $x \in \mathcal{M}$, the semi-supervised Nadaraya-Watson (NW) Fr\'echet regressor at $x$ is
	\begin{align}\label{Semi NW estimator}
		\hat{m}_{x}^{\textsl{snw}}=\underset{y \in \Omega}{\argmin} \ \frac{1}{n} \sum_{i=1}^{n} K_{h}\left(d_G(X_i, x)\right) d_{\Omega}^{2}\left(Y_i, y\right),
	\end{align}
	where $d_G(X_i, x)$ can be acquired by \eqref{graph distance}, $K$ is a smoothing kernel such as the  Epanechnikov kernel or Gaussian Kernel and $h$ is a bandwidth, with $K_{h}(\cdot)=h^{-p} K(\cdot / h)$. When $\Omega=\mathcal{R}$, we can verify that
	\begin{align*}
		\hat{m}^{\textsl{snw}}_x=\frac{ \sum_{i=1}^{n} K_{h}\left(d_G(X_i, x)\right) Y_i}{ \sum_{i=1}^{n} K_{h}\left(d_G(X_i, x)\right)} ,
	\end{align*}
	which has an explicit expression.
	
	Additionally, we define the semi-supervised $k$-nearest neighbor (kNN) Fr\'echet regressor at $x$ as
	\begin{align}\label{Semi kNN estimator}
		\hat{m}^{\textsl{skn}}_{x}=\underset{y \in \Omega}{\argmin} \sum_{i : X_i \in \textup{kNN}(x, d_G)} d_{\Omega}^2(Y_i, y),
	\end{align}
	where $\textup{kNN}(x, d_G)$ represents the set of $k$ nearest labeled neighbors to $x$ based on the graph distance $d_G$. Likewise, when $\Omega=\mathcal{R}$, it has the explicit expression
	\begin{align*}
		\hat{m}^{\textsl{knn}}_{x}=\frac{1}{k} \sum_{i : X_i \in \textup{kNN}(x, d_G)}  Y_i,
	\end{align*}
	which is just the Euclidean semi-supervised method proposed by \cite{moscovich2017minimax}.

	In essence, the workflow begins by constructing a graph that serves as a discrete representation of the underlying manifold, incorporating all available $X$ samples. The graph is then employed for the subsequent unsupervised metric learning. Following that, the Euclidean distances in the supervised Fr\'echet regression methods are replaced with estimated geodesic distances, resulting in a semi-supervised adaptation of Fr\'echet regression.  The abundant unlabeled data help accurately estimate the geodesic distances between the sample points on a low-dimensional manifold, as if the geodesic distances were explicitly known.  From this perspective, our semi-supervised Fr\'echet regression is a regression problem from a manifold space to a general metric space. 
	
	It is worth noting that besides the Fr\'echet regression itself, our implementation of inductive semi-supervised regression differs from \cite{moscovich2017minimax}. In their approach, the response values of all vertices in the original graph are first predicted by transductive learning. When handling a new coming point, the predicted response value of the vertex closest to the point (in terms of Euclidean distance) is selected as its prediction. Although this realization of inductive learning simplifies computational complexity, it introduces additional bias and raises challenges for theoretical analysis due to the potential mismatch between nearest neighbors under Euclidean distances and geodesic distances. In contrast, our approach completely avoids this issue by adding the new coming point to the original graph. Because of this, different points will correspond to different graphs, and the classical graph approximation theory \citep{tenenbaum2000global} used in \cite{moscovich2017minimax} will no longer be applicable to the inductive learning strategy here. We will tackle this theoretical difficulty from the perspective of the Hausdorff distance between sets of points. And our theoretical guarantees only require that the support set of $X$ is a compact manifold without imposing any assumptions on its density function, except for the implementation of the kNN algorithm on the manifolds of dimension $q < 3$.

	\section{Nonasymptotic analysis}\label{Convergence rate}
	To characterize the complexity of the response variable space, we employ Talagrand's $\gamma_2$ measure \citep{talagrand2014upper}. Below is the definition given in \cite{schotz2022nonparametric}. For the sake of completeness, we restate it here.
	\begin{definition}
		(i) Given a set $\mathcal{B}$, an admissible sequence is an increasing sequence $\left(\mathcal{A}_k\right)_{k \in \mathcal{N}_0}$ of partitions of $\mathcal{B}$ such that $\mathcal{A}_0=\{\mathcal{B}\}$ and the cardinality of $\mathcal{A}_k$ is bounded as $\# \mathcal{A}_k \leq 2^{2^k}$ for $k \geq 1$.
		By an increasing sequence of partitions, we mean that every set of $\mathcal{A}_{k+1}$ is contained in a set of $\mathcal{A}_k$. We denote by $A_k(\omega)$ the unique element of $\mathcal{A}_k$ which contains $\omega \in \mathcal{B}$.
		(ii) Let $(\mathcal{B}, d)$ be a pseudo-metric space, i.e., $d$ is symmetric, fulfills the triangle inequality, and $d(\omega, \omega)=0$ for all $\omega \in \mathcal{B}$. Define
		$$
		\gamma_2(\mathcal{B}, d):=\inf \sup _{\omega \in \mathcal{B}} \sum_{k=0}^{\infty} 2^{\frac{k}{2}} \operatorname{diam}\left(A_k(\omega), d\right)
		$$
		where the infimum is taken over all admissible sequences in $\mathcal{B}$, $\operatorname{diam}\left(A_k(\omega), d\right)$ represents the diameter of $A_k(\omega)$ with respect to the metric $d$.    
	\end{definition}
	
	Now we make the following assumptions:
	\begin{itemize}[leftmargin=0.8cm]
		\item[(A1)]
		\textsc{VarIneq}:
		There is $C_{\textsl{Vlo}}\in[1,\infty)$ such that $C_{\textsl{Vlo}}^{-1}\,d_{\Omega}^2(\omega, m_x) \leq M(x, \omega) - M(x, m_x)$ for all $\omega \in \Omega$ and $x \in \mathcal{M}$.
		\item[(A2)]
		\textsc{Entropy}: 
		There are $C_{\textsl{Ent}} \in [1,\infty)$ and $\alpha \in [1,2)$ such that 
		\begin{equation*}
			\gamma_2(\mathcal{B}, d_{\Omega}) \leq C_{\textsl{Ent}}\max\big(\diam(\mathcal{B}, d_{\Omega}), \diam(\mathcal{B}, d_{\Omega})^\alpha\big)
		\end{equation*} 
		for all $\mc B \subset \Omega$.
		\item[(A3)]
		\textsc{Moment}:
		There are $\kappa > \frac{2}{2-\alpha}$ and $C_{\textsl{Mom}} \in [1,\infty)$ such that $\E\big\{d_{\Omega}^\kappa(Y, m_x) \big| X=x\big\}^{\frac{1}{\kappa}} \leq C_{\textsl{Mom}}$ for all $x \in \mathcal{M}$.
		\item[(A4)] 
		\textsc{Kernel}:
		$K$ is a positive kernel on $[0, 1]$ highest at $0$, decreasing on $[0, 1]$ and $0$ outside $[0, 1]$.
		\item[(A5)] 
		\textsc{HölderSmoothDensity}:
		The function $\mathcal{M} \to \Omega,\, x \mapsto m_x$ is continuous. Let $C_{\textsl{Len}} \in [1,\infty)$ such that $\sup_{x_1, x_2\in\mathcal{M}} d_{\Omega}(m_{x_1}, m_{x_2}) \leq C_{\textsl{Len}}$. 
		Let $\mu$ be a probability measure on $\Omega$. Let $C_{\textsl{Int}} \in [1,\infty)$ such that $\int d_{\Omega}^2(y, m_0) \mu (\dl y) \leq C_{\textsl{Int}}$. 
		Let $y \to \rho(y|x)$ be the $\mu$-density of $Y$ conditional on $X=x$.
		Let $\beta \in (0,1]$. For $\mu$-almost all $y \in \Omega$, there is $L(y) \geq 0$ such that 
		$$\left|\rho(y \mid x_1)- \rho(y \mid x_2)\right| \leq  L(y) \cdot d_{\mathcal{M}}^\beta(x_1, x_2).$$ Furthermore, there is a constant $C_{\textsl{SmD}} > 0$, $\int L(y)^2  \mu(\dl y) \leq C_{\textsl{SmD}}^2$.
		\item[(A6)]
		\textsc{BiasMoment}: 
		Define $H(\omega_1, \omega_2) = \big\{\int  \big(d_{\Omega}(y, \omega_1) + d_{\Omega}(y, \omega_2)\big)^2 \mu (\dl y)\big\}^\frac{1}{2}$. There is $C_{\textsl{Bom}}\in[1,\infty)$ such that $\E \big\{H(\hat m_x, m_x)^\kappa \big| \{X_i\}_{i=1}^N\big\}^{\frac{1}{\kappa}} \leq C_{\textsl{Bom}}$ for all $x\in \mathcal{M}$.
		\item[(A7)]  \textsc{Manifold}: The intrinsic and ambient topologies coincide on $\mathcal{M}$ and the shortest paths on $\mathcal{M}$ have curvature bounded by $\beta$.
	\end{itemize}
	The above assumptions are the standard assumptions initially introduced by \cite{schotz2022nonparametric} to establish the convergence rate of local Fr\'echet regression in expectation. These assumptions have been relaxed as far as possible and are justified by demonstrating their fulfillment on hyperspheres. Therefore, we continue to adopt these very assumptions in our analysis to bolster the persuasiveness of our arguments. 
	The assumption (A1) is a common condition that controls the convergence rate of M-estimators. And it always holds true when $\Omega$ is a Hadamard space (complete geodesic spaces with curvature upper bounded by $0$), which includes Euclidean space $\mathcal{R}^p$ as a special case. 
	For the assumption (A2), it is satisfied with $\alpha=1$ for both the Euclidean space and any bounded metric space (metric space with bounded diameter). The bounded Fr\'echet moment condition (A3) is a natural generalization of the moment requirement in the Euclidean case.   
	The assumption (A4) is a classical requirement for local kernel regression, which is different from that in \citep{schotz2022nonparametric} because here we consider the NW regression with random designs. 
	In order to restrict the function space and then achieve the nonparametric rate of convergence in Fr\'echet regression, the assumption (A5) requires the smoothness of conditional Fr\'echet mean by the smoothness of conditional density $\rho(y|x)$. This is imperative since the conditional Fr\'echet mean does not have an explicit expression. Specifically, when $\Omega=\mathcal{R}$, the H\"older continuity of the conditional density readily implies the H\"older continuity of the regression function. 
	Lastly, the assumption (A6) is a general condition that can be fulfilled by bounded metric space and Hadamard space. For more detailed comments on these assumptions, please refer to Remarks 1 and 2 in \cite{schotz2022nonparametric}. The assumption (A7) is some loose requirements on the geometrical properties of the manifold $\mathcal{M}$, which is a summary of Properties 3.4 and 3.11 from \cite{arias2019unconstrained}. The assumption (A7) is made to ensure that the geodesic distance can be controlled by the graph distance with an appropriate graph radius $r$.\\
	
	\begin{theorem}[Semi-supervised Fr\'echet regression]\label{semi-supervised theorem}
		Assume $X$ is supported on a $q$-dimensional submanifold $\mathcal{M}$ embedded in $\mathcal{R}^p$ $(q \leq p)$. Moreover, assume (A1)--(A7). Given  $\epsilon >0, \lambda \in (0,1)$ and $\tau, c_0 > 0$ defined in Lemma 6 of supplementary materials. Then for enough large $m$ satisfying $m^{(\lambda-1)/q} \leq \epsilon$ and $4 \epsilon \leq r \leq \min\big(\tau,  1/(3\beta)\big)$, let $\delta_1=\frac{1}{1+c_0 r^2}$ and $\delta_2=1+4 \varepsilon / r$, the following results holds based on the $r$-graph constructed by $\{X_i\}_{i=1}^N$
		\begin{itemize}
			\item For the NW semi-supervised Fr\'echet regression,
			\begin{align*}
				\int\E\left\{d_{\Omega}^2(\hat m_{x}^{\textsl{snw}}, m_{x})\right\} \nu_X(\dl x)
				\leq &
				c_{\alpha,\kappa} \br{C_{\textsl{Vlo}} C_{\textsl{SmD}} C_{\textsl{Bom}}}^{\frac{2}{2-\alpha}}h^{2\beta}/ \delta_1^{2\beta}+ 
				\\&\hphantom{\leq}\ \,c_{\alpha,\kappa} \br{C_{\textsl{Vlo}}C_{\textsl{Mom}}C_{\textsl{Len}}}^{\frac{2}{2-\alpha}} e^{-1} \delta_2^q/(nh^q)+
				\\&\hphantom{\leq}\ \, c_{\alpha,\kappa}
				\br{ C_{\textsl{Vlo}} C_{\textsl{Mom}} C_{\textsl{Ent}} }^{\frac{2}{2-\alpha}} \delta_2^q/(nh^q)+
				\\&\hphantom{\leq}\ \, c_{\alpha,\kappa} \br{C_{\textsl{Vlo}} C_{\textsl{SmD}} C_{\ms{Bom}}C_{\textsl{Mom}}C_{\textsl{Len}}C_{\textsl{Ent}}}^{\frac{2}{2-\alpha}} \diam(\mathcal{M})^{2\beta} \exp \left(-\theta m^\lambda\right);
			\end{align*}
			Taking $h \asymp n^{-\frac{1}{2\beta+q}}$, then we get
			\begin{align*}
				\int\E\left\{d_{\Omega}^2(\hat m_{x}^{\textsl{snw}}, m_{x})\right\} \nu_X(\dl x) \leq c \cdot n^{-\frac{2\beta}{2\beta+q}}+c' \cdot \exp \left(-\theta m^\lambda\right).
			\end{align*}
			\item For the kNN semi-supervised Fr\'echet regression, if $q \geq 3$,
			\begin{align*}
				\int\E\left\{d_{\Omega}^2(\hat m_x^{\textsl{skn}}, m_x)\right\} \nu_X(\dl x)
				\leq &
				c_{\alpha,\kappa} \br{C_{\textsl{Vlo}} C_{\textsl{SmD}} C_{\textsl{Bom}}}^{\frac{2}{2-\alpha}} \left(\frac{\delta_2}{\delta_1}\right)^{4\beta} \frac{c}{\left\lfloor n/k \right\rfloor^{2\beta / q}} +
				\\&\hphantom{\leq}\ \, c_{\alpha,\kappa}
				\br{ C_{\textsl{Vlo}} C_{\textsl{Mom}} C_{\textsl{Ent}} }^{\frac{2}{2-\alpha}} k^{-1} +
				\\&\hphantom{\leq}\ \, c_{\alpha,\kappa} \br{C_{\textsl{Vlo}} C_{\textsl{SmD}} C_{\ms{Bom}}}^{\frac{2}{2-\alpha}} \diam(\mathcal{M})^{2\beta} \exp \left(-\theta m^\lambda\right);
			\end{align*}
			if $p < 3$, the above bound still holds provided there is one more assumption that there exist $\xi_0>0$, a nonnegative function $g$ such that for all $x \in \mathcal{M}$, and $0<\xi \leq \xi_0$,
			\begin{align*}
				\nu_X\big(S(x, d_{\mathcal{M}}, \xi)\big)>g(x) \xi^q \quad \text{and} \quad  \int \frac{1}{g(x)^{2 / q}} \nu_X(\dl x)<\infty.
			\end{align*}
			Taking $k \asymp n^{\frac{2\beta}{2\beta+q}}$, then we get
			\begin{align*}
				\E\left\{d_{\Omega}^2(\hat m_x^{\textsl{skn}}, m_x)\right\} \leq c \cdot n^{-\frac{2\beta}{2\beta+q}}+c' \cdot \exp \left(-\theta m^\lambda\right).
			\end{align*}
		\end{itemize}
	\end{theorem}
	
	We are mainly interested in the convergence rate and do not explore the best universal constants here. The exponential term of the upper error bound above is the manifold approximation error which is negligible compared to the first term, when the size of unlabeled data $m$ is large enough. Therefore
	Theorem~\ref{semi-supervised theorem} reveals that the convergence rate of two semi-supervised Fr\'echet regression methods adapts to the intrinsic dimension $q$ of the submanifold $\mathcal{M}$, rather than the ambient dimension $p$, even if the number $n$ of labeled samples is finite.  In particular, the rate $n^{-\frac{2\beta}{2\beta+q}}$ is minimax optimal \citep{stone1982optimal} for nonparametric estimation of $\beta$- H\"older continuous regression functions on $\mathcal{R}^q$.
	
	\begin{remark}
		Like \cite{schotz2022nonparametric}, we focus on two important metric spaces: bounded metric spaces and Hadamard space.
		\begin{itemize}
			\item[(1)] When $\Omega$ is a bounded metric space, the assumption (A2) holds true with $\alpha=1$; the assumption (A3) holds true with $C_{\textsl{Mom}}=4 \diam(\Omega, d_{\Omega})$; the assumption(A6) holds true with $C_{\textsl{Bom}}=2 \diam(\Omega, d_{\Omega})$.
			\item[(2)] When $\Omega$ is a Hadamard space, the assumption (A1) holds true with $C_{\textsl{Vlo}}=1$; the assumption (A6) holds true with $C_{\textsl{Bom}}=c_{\kappa}C_{\textsl{Mom}}C_{\textsl{Len}}C_{\textsl{Int}}$ by Proposition 3 of \cite{schotz2022nonparametric}.
		\end{itemize}
	\end{remark}
	
	As a byproduct of the above analysis, we consider to establish the nonasymptotic convergence rate of excess risk for supervised NW Fr\'echet regression
	\begin{align}\label{NW estimator}
		\hat{m}_{x}^{\textsl{nw}}=\underset{\omega \in \Omega}{\argmin} \frac{1}{n} \sum_{i=1}^{n} K_{h}\left(\|X_{i}-x\|_2\right) d_{\Omega}^{2}\left(Y_{i}, \omega\right).
	\end{align}
	and supervised kNN Fr\'echet regression
	\begin{align}\label{kNN estimator}
		\hat{m}^{\textsl{knn}}_{x}=\underset{\omega \in \Omega}{\argmin} \sum_{i: X_i \in \textup{kNN}(x, \|\cdot\|_2)} d_{\Omega}^2(Y_i, \omega).
	\end{align}
	For the supervised scenario, \eqref{NW estimator} and \eqref{kNN estimator} are local methods developed based on the Euclidean distance. The adaptions of the assumptions (A5) and (A6) are needed as
	
	\begin{itemize}
		\item[(B5)] 
		\textsc{HölderSmoothDensity}:
		The function $[0,1]^p \to \Omega,\, x \mapsto m_x$ is continuous. Let $C_{\textsl{Len}} \in [1,\infty)$ such that $\sup_{x_1, x_2\in[0,1]^p} d_{\Omega}(m_{x_1}, m_{x_2}) \leq C_{\textsl{Len}}$. 
		Let $\mu$ be a probability measure on $\Omega$. Let $C_{\textsl{Int}} \in [1,\infty)$ such that $\int d_{\Omega}^2(y, m_0) \mu (\dl y) \leq C_{\textsl{Int}}$. 
		Let $y \to \rho(y|x)$ be the $\mu$-density of $Y$ conditional on $X=x$.
		Let $\beta \in (0,1]$. For $\mu$-almost all $y \in \Omega$, there is $L(y) \geq 0$ such that 
		$$\big|\rho(y | x_1)- \rho(y | x_2)\big| \leq  L(y) \cdot \|x_1 - x_2\|^{\beta}.$$ Furthermore, there is a constant $C_{\textsl{SmD}} > 0$, $\int L(y)^2  \mu(\dl y) \leq C_{\textsl{SmD}}^2$.
		\item[(B6)]
		\textsc{BiasMoment}: 
		Define $H(\omega_1, \omega_2) = \big\{\int  \big(d_{\Omega}(y, \omega_1) + d_{\Omega}(y, \omega_2)\big)^2 \mu (\dl y)\big\}^\frac{1}{2}$. There is $C_{\textsl{Bom}}\in[1,\infty)$ such that $\E \big\{H(\hat m_x, m_x)^\kappa \big| \{X_i\}_{i=1}^n\big\}^{\frac{1}{\kappa}} \leq C_{\textsl{Bom}}$ for all $x\in[0,1]^p$.
	\end{itemize} 
	
	\begin{theorem}[Supervised Fr\'echet regression]\label{supervised theorem}
		Assume (B5)–(B6) and (A1)–(A4) with $x \in \mathcal{M}$
		replaced by $x \in [0, 1]^p$. Then it holds that
		\begin{itemize}
			\item For the NW Fr\'echet regression,
			\begin{align*}
				\int\E\left\{d_{\Omega}^2(\hat m_x^{\textsl{nw}}, m_x)\right\} \nu_X(\dl x)
				\leq &
				c_{\alpha,\kappa} \br{C_{\textsl{Vlo}} C_{\textsl{SmD}} C_{\textsl{Bom}}}^{\frac{2}{2-\alpha}} h^{2\beta} + 
				\\&\hphantom{\leq}\ \,c_{\alpha,\kappa} \br{C_{\textsl{Vlo}}C_{\textsl{Mom}}C_{\textsl{Len}}}^{\frac{2}{2-\alpha}} e^{-1} (nh^p)^{-1}+
				\\&\hphantom{\leq}\ \, c_{\alpha,\kappa}
				\br{ C_{\textsl{Vlo}} C_{\textsl{Mom}} C_{\textsl{Ent}} }^{\frac{2}{2-\alpha}} (nh^p)^{-1}.
			\end{align*}
			Taking $h \asymp n^{-\frac{1}{2\beta+p}}$, then we get
			\begin{align*}
				\int\E\left\{d_{\Omega}^2(\hat m_x^{\textsl{nw}}, m_x)\right\} \nu_X(\dl x) \leq c \cdot n^{-\frac{2\beta}{2\beta+p}}.
			\end{align*}
			\item For the kNN Fr\'echet regression, if $p \geq 3$,
			\begin{align*}
				\int\E\left\{d_{\Omega}^2(\hat m_x^{\textsl{kn}}, m_x)\right\} \nu_X(\dl x)
				\leq &
				c_{\alpha,\kappa} \br{C_{\textsl{Vlo}} C_{\textsl{SmD}} C_{\textsl{Bom}}}^{\frac{2}{2-\alpha}} \frac{c}{\left\lfloor n/k \right\rfloor^{2\beta / p}}
				+
				\\&\hphantom{\leq}\ \, c_{\alpha,\kappa}
				\br{ C_{\textsl{Vlo}} C_{\textsl{Mom}} C_{\textsl{Ent}} }^{\frac{2}{2-\alpha}} k^{-1};
			\end{align*}
			if $p < 3$, the above bound still holds provided there is one more assumption that there exist $\xi_0>0$, a nonnegative function $g$ such that for all $x \in [0,1]^p$, and $0<\xi \leq \xi_0$,
			\begin{align*}
				\nu_X\big(S(x, \|\cdot\|_2, \xi)\big)>g(x) \xi^p \quad \text{and} \quad  \int \frac{1}{g(x)^{2 / p}} \nu_X(\dl x)<\infty.
			\end{align*}
			Taking $k \asymp n^{\frac{2\beta}{2\beta+p}}$, then we get
			\begin{align*}
				\int\E\left\{d_{\Omega}^2(\hat m_x^{\textsl{kn}}, m_x)\right\} \nu_X(\dl x) \leq c \cdot n^{-\frac{2\beta}{2\beta+p}}.
			\end{align*} 
		\end{itemize}
	\end{theorem}
	The above theorem shows that both supervised Fr\'echet regression methods achieve the same convergence rate. The bound here is regarding the convergence in expectation and is nonasymptotic, unlike the results of \cite{petersen2019frechet}. In comparison with Theorem~\ref{supervised theorem}, it signifies that the intrinsic low-dimensional manifold structure of the feature space is a fundamental factor enabling the significant improvement of convergence rate in semi-supervised methods. What deserves attention here is the work of \cite{bickel2007local} when considering local polynomial Euclidean regression on unknown manifolds. It reveals ``naive'' local polynomial regression can adapt to local smooth lower dimensional structure in the sense that its asymptotic convergence rate is determined by the intrinsic $q$ rather than ambient dimension $p$. Leaving aside the difference that we are studying Fr\'echet regression, their work seems to indicate that our semi-supervised treatment is redundant. This is not actually the case. Their theoretical results rely on a situation where there are enough labeled samples, leading to, locally, the geodesic distance being roughly proportionate to the Euclidean distance. Under this premise,  they can unwittingly take advantage of low-dimensional structure without manifold estimation.  In contrast, the result of Theorem~\ref{supervised theorem} is nonasymptotic, allowing the number of labeled samples $n$ to be any finite positive integer. Consequently, we necessitate additional unlabeled data to assist in estimating the manifold structure, and then utilize the acquired structural information (geodesic distance) to achieve effective supervised learning on limited labeled samples. Therefore, our semi-supervised processing retains significance. Of course, according to \cite{bickel2007local}, semi-supervised learning becomes unnecessary when the number of labeled samples is sufficiently large.

	\section{Simulation}\label{Simulation}
	To demonstrate the efficacy of our proposed semi-supervised methods, we apply kNN Fr\'echet regression \eqref{kNN estimator}, semi-supervised kNN Fr\'echet regressionc \eqref{Semi kNN estimator}, NW Fr\'echet regression \eqref{NW estimator} and semi-supervised NW Fr\'echet regression \eqref{Semi NW estimator}  on a series of simulations.
	Based on the labeled data $\{(X_i, Y_i)\}_{i=1}^n$ with size $n$ and unlabeled data $\{X_i\}_{i=n+1}^{N}$ with size $m=N-n$, the four methods are trained to predict responses for another independent testing data $\{X_i\}_{i=N+1}^{N+1000}$ with size $1000$. But the unlabeled data are not used for the two supervised learning methods, i.e., kNN Fr\'echet regression and NW Fr\'echet regression.  For any simulation setting, we evaluate the performance of each method by computing the average mean squared error (AMSE) over $100$ realizations. Specifically, for the $j$th Monte Carlo realization, $\hat{m}^j$ denotes the fitted Fr\'echet regression function based on the method $\hat{m}$ and the quality of prediction is measured by the mean squared error 
	$$
	\operatorname{MSE}_{j}(\hat{m})=  \frac{1}{1000} \sum_{i=N+1}^{N+1000} d_{\Omega}^2\big(\hat{m}^j_{X_i},m_{X_i}\big)
	$$
	based on the testing points. And the average mean squared error of the method $\hat{m}$ is calculated by
	$$
	\operatorname{AMSE}(\hat{m})=\frac{1}{100} \sum_{j=1}^{100} \operatorname{MSE}_{j}(\hat{m}).
	$$
	
	For NW Fr\'echet regression and semi-supervised NW Fr\'echet regression,  the kernel function $K$ is taken as the Epanechnikov kernel $K(u)=3 / 4  \cdot \big(1 -  u^2\big) \ind_{\{u \in [-1,1]\}}$. For two semi-supervised methods, there are no universally clear guidelines for selecting the connectivity radius $r$ of $r$-graph. Often, it can be selected via cross-validation within an appropriate range. Here $r$ is set as 1.2 times the maximal value of the minimal Euclidean distances of each point in $\{X_i\}_{i=1}^{N}$ to the other points, that is,
	$$ 
	r=1.2 \cdot \max_{1 \leq i \leq N}\big\{\min_{1 \leq j \neq i \leq N}\{\|X_i -  X_j\|_2\big\}.
	$$ 
	This empirical choice of $r$ behaves particularly well in all of the following simulations. Lastly, we use leave-one-out cross validation to select the hyperparameter $k$ among $\{1,2, \ldots, 10\}$ for two kNN methods and bandwidth $h$ among $\{5^{0.1} \cdot h_0, 5^{0.2} \cdot h_0, \ldots, 5^{1} \cdot h_0\}$ for two NW methods, where $h_0$ is the the median value of the minimal Euclidean distances (for NW Fr\'echet regression) or graph distances (for semi-supervised NW Fr\'echet regression) of each points in $\{X_i\}_{i=1}^{n}$ to the other points.
	
	\subsection{Fr\'echet regression on $\mathcal{R}^3$}\label{three dimension}
\begin{figure}[ht!]
		\centering
		\includegraphics[width=0.4\linewidth]{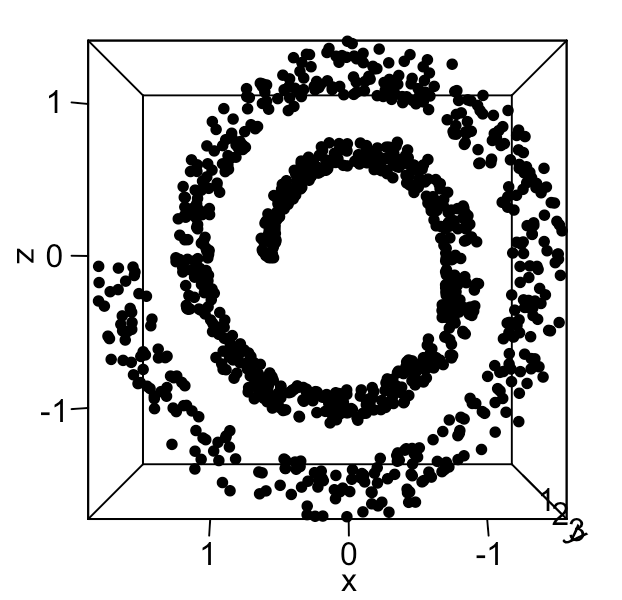}
		\caption{A typical Swiss roll in three-dimensional Euclidean space with $1000$ sample points.}
		\label{fig:swiss-roll}
	\end{figure}

	We start with simulated data where $X$ is distributed on a Swiss roll, an intrinsically 2-dimensional manifold in the 3-dimensional Euclidean space (as depicted in Figure~\ref{fig:swiss-roll}). $X=\left(X^{(1)}, X^{(2)}, X^{(3)}\right)^\T$ is generated by
	$$
	X^{(1)}=\theta \cos (\theta)/10, \quad X^{(2)}=4U^{(2)}, \quad X^{(3)}=\theta \sin (\theta)/10,
	$$
	with $\theta = 4\pi(U^{(1)}+1/2)$ and $U=\big(U^{(1)}, U^{(2)}\big)^{\T} \sim \operatorname{Unif}([0,1]^2)$. We first sample $N$ points uniformly from $[0,1]^2$, denoted as $\left\{U_i\right\}_{i=1}^N$, and then obtain the corresponding $N$ observations $\left\{X_i\right\}_{i=1}^N$ on $X$ by the above parameterization. It is clear that the sampling design on the $X$ is not uniform even if $U$ is sampled uniformly. The nonEuclidean response $\{Y_i\}_{i=1}^n$ depending on $\{X_i\}_{i=1}^n$ are generated by unobservable $\{U\}_{i=1}^n$. The detailed generation process will be described in Sections 4.1.1 and 4.1.2. The labeled data size is set to be $n=100$, and unlabeled data size $m$ is considered among $\{0, 100, 200, 500, 1000, 2000, 3000\}$.  
	
	\emph{4.1.1. Responses for symmetric positive-definite matrices.}
	Let $(\Omega,d_{\Omega})$  be the metric space  $\mathcal{S}_m^{+}$ of $m \times m$ symmetric positive-definite matrices endowed with metric $d_{\Omega}$. There are many options for metrics, this section focuses on the Log-Cholesky metric \citep{lin2019riemannian}. For a matrix $Y$, let $\lfloor Y\rfloor$ denote the strictly lower triangular matrix of $Y$,  $\mathbb{D}(Y)$ denote the diagonal part of $Y$ and $\|Y\|_{F}$ denote the Frobenius norm. It is well known that if $Y$ is a symmetric positive-definite matrix, there is a lower triangular matrix $P$ whose diagonal elements are all positive such that $PP^{T}=Y$.  This $P$ is called the Cholesky factor of $Y$, devoted as $\mathscr{L}(Y)$. For an $m \times m$ symmetric matrix $A$, $\exp (A)=I_{m}+\sum_{j=1}^{\infty} \frac{1}{j !} A^{j}$ is a symmetric positive-definite matrix. Conversely, for a symmetric positive-definite matrix $Y$, the matrix logarithmic map is $\log (Y)= A$ such that $\exp (A)=Y$. For two symmetric positive-definite matrices $y_{1}$ and $y_{2}$, the Log-Cholesky metric is defined by
	\begin{align*}
		d_{\Omega}(y_1, y_2)=d_{\mathcal{L}}(\mathscr{L}(y_1), \mathscr{L}(y_2)),
	\end{align*}
	where $d_{\mathcal{L}}(P_1,P_2)=\{\left\| \lfloor P_1 \rfloor-\lfloor P_2 \rfloor\right\|_{F}^{2}+\left\| \log \mathbb{D}(P_1)-\log \mathbb{D}(P_2)\right\|_{F}^{2}\}^{1 / 2}$.
	
	The response $Y$ is generated via symmetric matrix variate normal distribution \citep{zhang2021dimension}. Consider the simplest case, we say an $m \times m$ symmetric matrix $A \sim \mathcal{N}_{m m}(M ; \sigma^2)$ if $A=\sigma Z+M$ where $M$ is an $m \times m$ symmetric matrix  and  $Z$ is an  $m \times m$ symmetric random matrix with independent $\mathcal{N}(0, 1)$ diagonal elements and $ \mathcal{N}(0, 1/2)$ off-diagonal elements. The response $Y$ depending on $X$ is generated by unobservable $U$.  Specifically, we consider the following settings.
	
	Setting I:
	$$
	\log (Y) \sim \mathcal{N}_{m m}\left(D(X), \sigma^{2}\right)
	$$
	with 
	$D(X)= \left(\begin{array}{cc}
		1 & \rho(X) \\
		\rho(X) & 1
	\end{array}\right),  \rho(X)= \cos\left(4  \pi (\beta^{\T}U) \right), \beta=(0.75, 0.25).$
	
	Setting II:
	$$
	\log (Y) \sim \mathcal{N}_{m m}\left( D(X), \sigma^{2}\right)
	$$
	with $D(X)=\left(\begin{array}{ccc}
		1 & \rho_{1}(X) & \rho_{2}(X) \\
		\rho_{1}(X) & 1 & \rho_{1}(X) \\
		\rho_{2}(X) & \rho_{1}(X) & 1
	\end{array}\right), \rho_1(X)=0.8 \cdot \cos\left(4  \pi (\beta_1U) \right), \rho_2(X)= 0.4 \cdot \cos\left(4  \pi (\beta_2U) \right), \beta_1=(0.75, 0.25), \beta_2=(0.25, 0.75).$
	
\begin{figure}[ht!]
		\centering
		\subfloat[setting I under snr=2]{\includegraphics[width=0.5\linewidth]{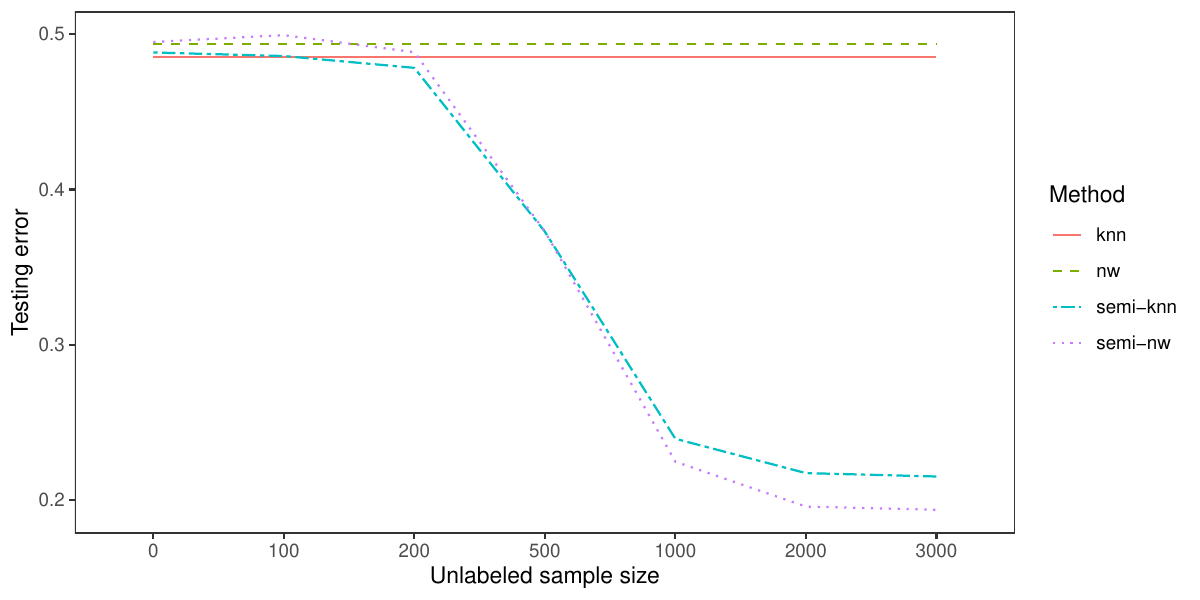}}
		\subfloat[setting I under snr=4]{\includegraphics[width=0.5\linewidth]{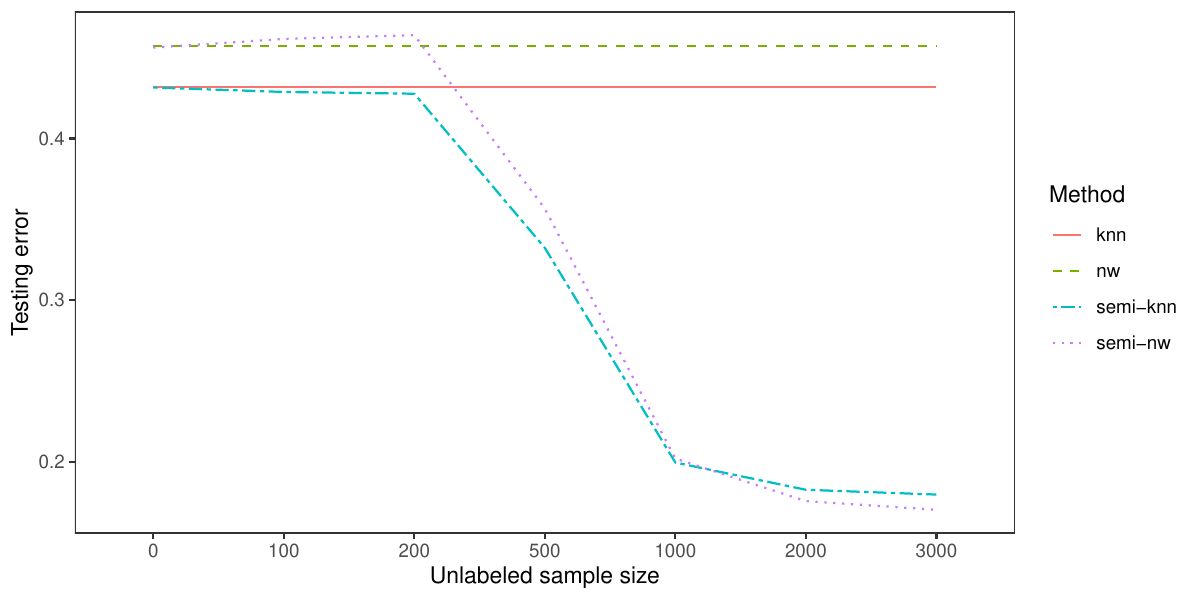}}
		\vspace{4mm}
		\subfloat[setting II under snr=2]{\includegraphics[width=0.5\linewidth]{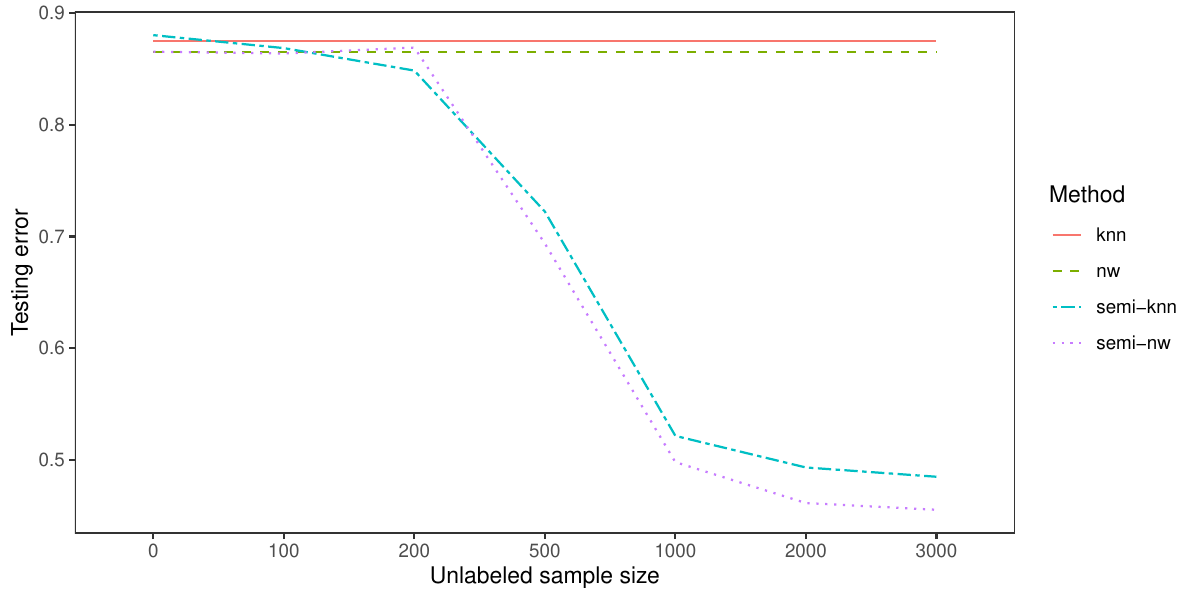}}
		\subfloat[setting II under snr=4]{\includegraphics[width=0.5\linewidth]{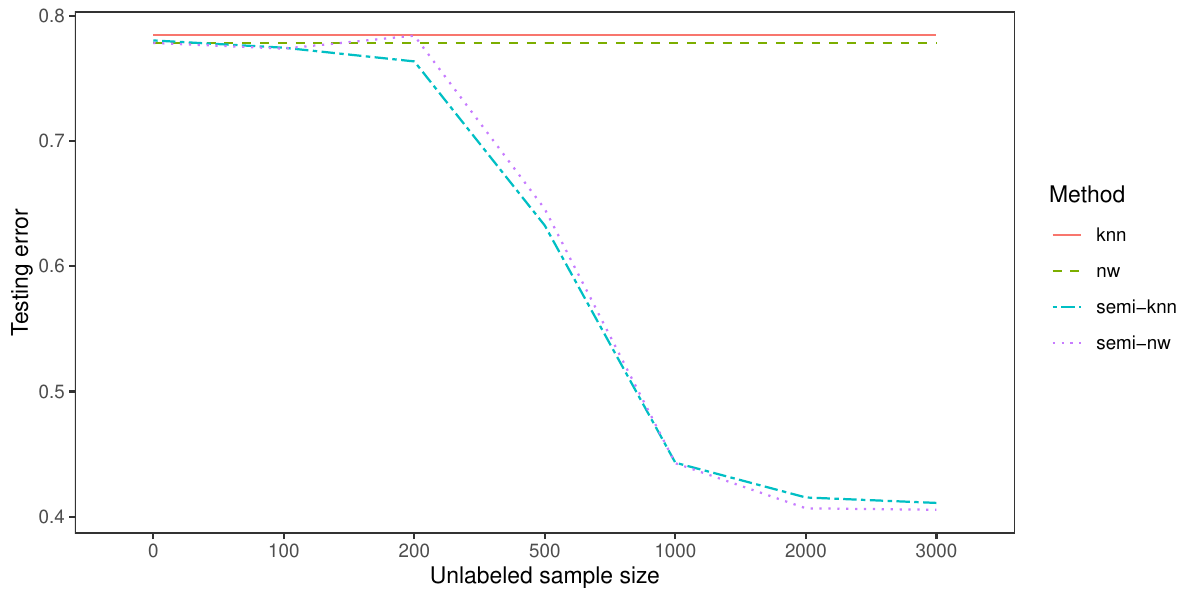}}
		\caption{AMSE of different methods for setting I, II with $100$ labeled points under snr$=2, 4$ when the dimension of $X$ is equal to $3$.}
		\label{fig: spd_three dimension}
	\end{figure}

	In these two settings, we choose the appropriate $\sigma$ to make the signal-to-noise ratio (snr) equal to $2$ (high noise level) or $4$ (low noise level).
	Figure~\ref{fig: spd_three dimension} depicts the results obtained from four combinations of settings and noise levels.  Across all scenarios, our focus remains on the impact of varying sizes of unlabeled data on the semi-supervised approaches. Regarding the supervised methods, their performance remains unchanged since they solely rely on the labeled data. Furthermore, a consistent observation is that when the size of unlabeled data is limited (including scenarios without unlabeled data), the performance of the semi-supervised methods resembles that of the supervised methods, and occasionally even shows a slight degradation. This phenomenon can be rationalized by acknowledging that there is a substantial error associated with approximating the geodesic distances on a low-dimensional manifold by the graph distances when feature samples are not enough. Essentially, the semi-supervised methods struggle to accurately capture the structure of the low-dimensional manifold with small $m$. This limitation inhibits the utilization of local information in terms of geodesic distance to enhance the precision of response value predictions. Nevertheless, as the size of unlabeled data increases to a certain extent, the advantages of the semi-supervised methods become pronounced. Their AMSE experiences a substantial reduction as the number of unlabeled samples increases, eventually exhibiting a trend toward stabilization at a low error level. Notably, when the size of unlabeled data reaches $3000$, it becomes evident that the two semi-supervised methods surpass their two supervised counterparts by a substantial margin. Moreover, we find that the prediction performance of the semi-supervised methods is constantly approaching that of the supervised methods that use $\{(U_i, Y_i)\}_{i=1}^n$ to make predictions, provided that $\{U_i\}_{i=1}^n$ can be observed. At this juncture, the prediction error is primarily attributed to the presence of noise.
	
	Overall, as the noise level shifts from low (SNR=2) to high (SNR=4), or as the response variable changes from a $2 \times 2$ matrix to a $3 \times 3$ matrix, the prediction task becomes more challenging for all the methods, which can be reflected by the larger AMSE. When the size of unlabeled data is large enough, the semi-supervised NW Fr\'echet regression slightly outperforms the semi-supervised kNN Fr\'echet regression.
	
	\emph{4.1.2. Responses for spherical data.}
	Now we consider another type of response. Let $(\Omega,d_{\Omega})$ be the metric space $\mathbb{S}^2$ of sphere data endowed with the geodesic distance $d_{\Omega}$. For any two points $y_1, y_2 \in \mathbb{S}^2$, the geodesic distance is defined by
	\begin{align*}
		d_{\Omega} \left(y_{1}, y_{2}\right)=\arccos \left(y_{1}^{\T}y_{2}\right).
	\end{align*}
	And $Y$ is generated by the following two settings.
	
	Setting III:
	Let the Fr\'echet regression function be
	\begin{align*}
		m_X=\big(&\{1-(\beta_{1}^{\T}U)^{2}\}^{1 / 2} \cos (\pi (\beta_{2}^{\T}U)), \{1-(\beta_{1}^{\T}U)^{2}\}^{1 / 2} \sin (\pi (\beta_{2}^{\T}U)), \beta_{1}^{\T}U\big)^{\T},
	\end{align*}
	where $\beta_1=(1, 0),  \beta_2=(0, 1)$. We generate binary Normal noise $\varepsilon_{i}$ on the tangent space $T_{m_{X_i}} \mathbb{S}^2$, then map $\varepsilon_{i}$ back to $\mathbb{S}^2$ by Riemannian exponential map to get $Y_i$. Specifically, we first independently generate $\delta_{i 1}, \delta_{i 2} \stackrel{i i d}{\sim} \mathcal{N}\left(0,0.2^{2}\right)$, then let $\varepsilon_{i}=\delta_{i 1} v_{1}+\delta_{i 2} v_{2}$, where $\left\{v_{1}, v_{2}\right\}$ forms an orthogonal basis of tangent space $T_{m_{X_i}} \mathbb{S}^2$. Then $Y_{i}$ can be generated by
	$$
	Y_{i}=\operatorname{Exp}_{m_{X_{i}}}\left(\varepsilon_{i}\right)=\cos \left(\left\|\varepsilon_{i}\right\|_2\right) m_{X_{i}}+\sin \left(\left\|\varepsilon_{i}\right\|_2\right) \frac{\varepsilon_{i}}{\left\|\varepsilon_{i}\right\|_2}.
	$$
	
	Setting IV:
	Consider the following model
	\begin{align*}
		Y_{i}=\big(&\sin (\beta_{1}^{\T}U_i+\varepsilon_{i 1}) \sin (\beta_{2}^{\T}U_i+\varepsilon_{i 2}), \sin (\beta_{1}^{\T}U_i+\varepsilon_{i 1}) \cos (\beta_{2}^{\T}U_i+\varepsilon_{i 2}), \abs{\cos (\beta_{1}^{\T}U_i+\varepsilon_{i 1})}\big)^{\T},
	\end{align*}
	where the random noise $\varepsilon_{i 1}, \varepsilon_{i 2} \stackrel{i i d}{\sim} \mathcal{N}\left(0,0.2^{2}\right)$ are generated independently. The choice of $\beta_1, \beta_2$ is the same as setting III.
	
\begin{figure}[ht!]
		\centering 
		\subfloat[setting III]{\includegraphics[width=0.5\linewidth]{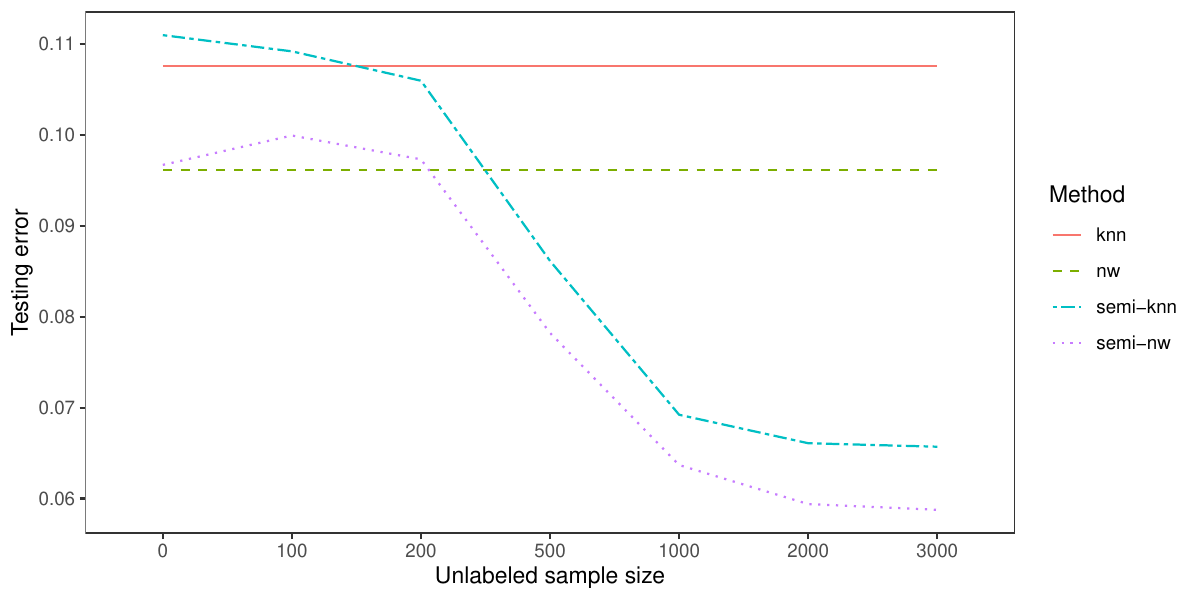}}
		\subfloat[setting IV]{\includegraphics[width=0.5\linewidth]{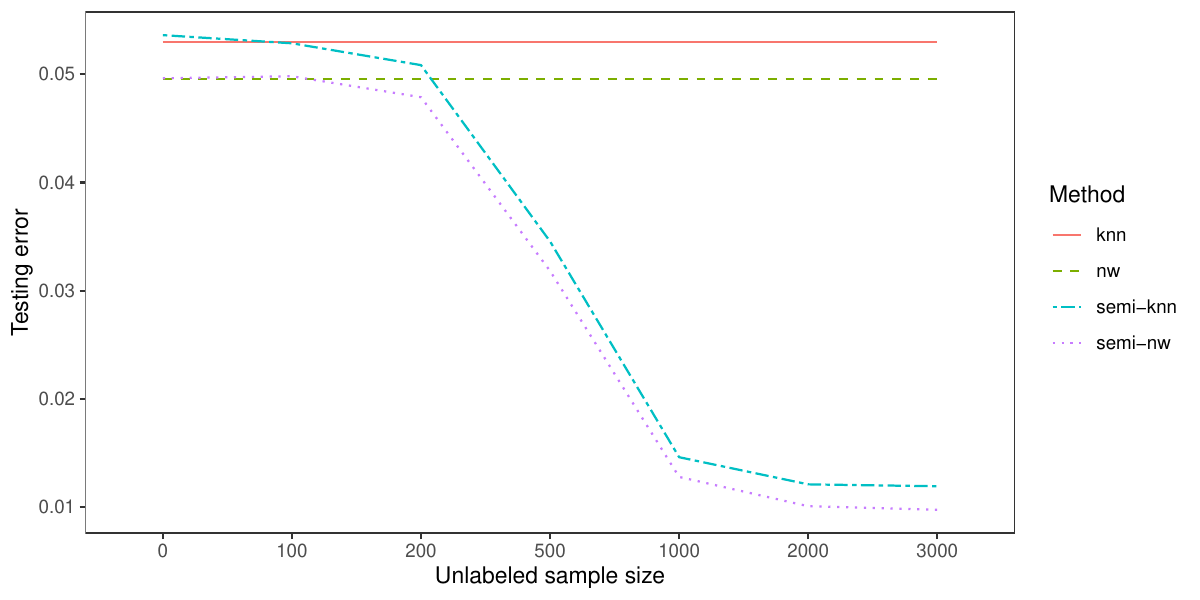}}
		\caption{AMSE of different methods for setting III, IV with $100$ labeled points when the dimension of $X$ is equal to $3$.}
		\label{fig:sph_theree dimension}
	\end{figure}

	The results are recorded in Figure~\ref{fig:sph_theree dimension}. The performance of all methods follows similar patterns as the simulations presented in the prior section. This observation underscores the general applicability of our semi-supervised approaches when addressing non-Euclidean problems, reinforcing the idea that leveraging low-dimensional manifold structures can significantly enhance prediction accuracy. Moreover, it can be found that in both settings, the NW methods consistently behave better than the KNN methods regardless of whether it is in the supervised or semi-supervised scenarios.

\begin{figure}[ht!]
		\centering
		\subfloat[setting I under snr=4]{\includegraphics[width=0.5\linewidth]{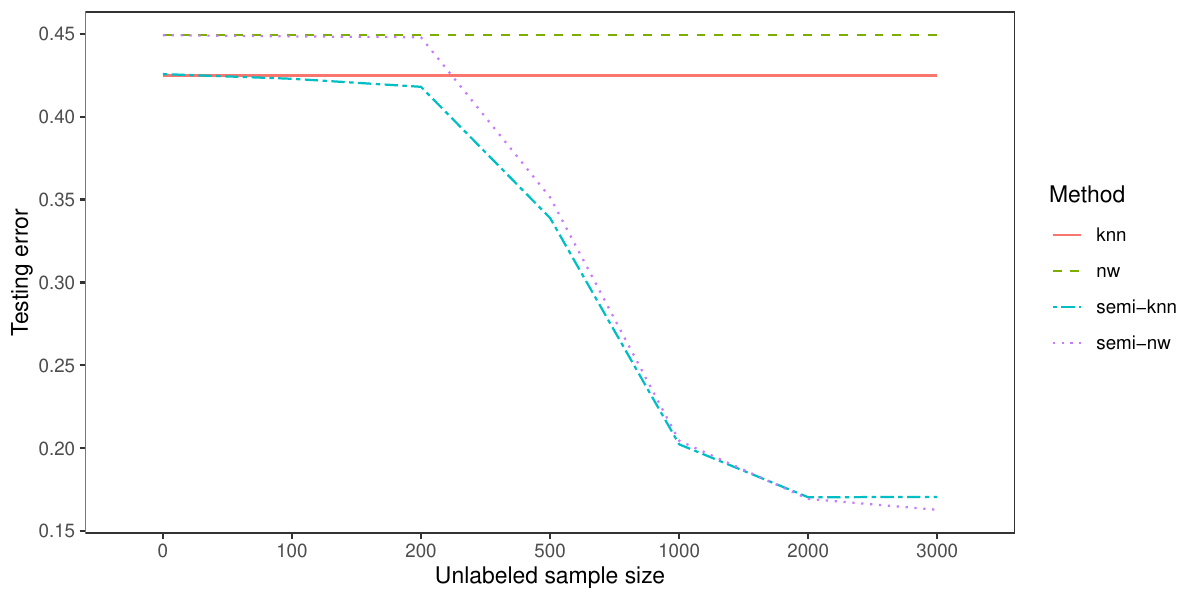}}
		\subfloat[setting III]{\includegraphics[width=0.5\linewidth]{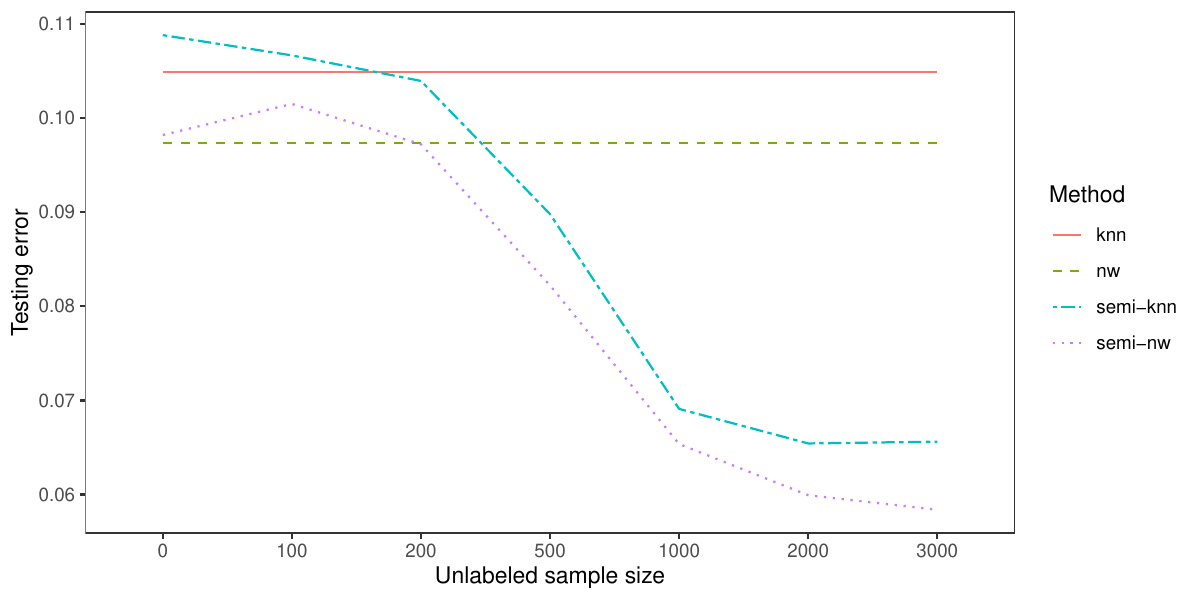}}
		\caption{AMSE of different methods for setting I, III with $100$ labeled points when the dimension of $X$ is equal to $3$ and the components of $U$ are correlated.}
		\label{fig:sph_three dimension_cor}
	\end{figure}

	Since the components of $U$ are independent in all previous simulations, we here sample $U$ from the truncated multivariate normal distribution supported on the compact set $[0,1]^2$ with mean vector $(0.5, 0.5)^\T$ and covariance matrix whose $(i, j)$ entry is $0.5^{|i-j|}$. We only consider the setting I under snr=4 and setting III. The respective outcomes are illustrated in Figure ~\ref{fig:sph_three dimension_cor} and similar to the case where $U$ comes from a uniform distribution.
	
	\subsection{Fr\'echet regression on $\mathcal{R}^6$}\label{six dimension}
	Let us consider another case where $X$ lies on a two-dimensional manifold embedded in $\mathcal{R}^6$. $X=\left(X^{(1)}, X^{(2)}, X^{(3)}, X^{(4)}, X^{(5)}, X^{(6)}\right)^\T$ can be parameterized by
	\begin{gather*}
		X^{(1)}=\theta_1 \cos (\theta_1)/10, \quad X^{(2)}=U^{(2)}, \quad X^{(3)}=\theta_1 \sin (\theta_1)/10,\\
		X^{(4)}=\theta_2 \cos (\theta_2)/10, \quad X^{(5)}=-U^{(1)}, \quad X^{(6)}=\theta_2 \sin (\theta_2)/10.
	\end{gather*}
	where $\theta_i = 4\pi(U^{(i)}+1/2), i=1,2$, and $U=\big(U^{(1)}, U^{(2)}\big)^{\T} \sim \operatorname{Unif}([0,1]^2)$. The data generation process parallels that outlined in Section~\ref{three dimension}, with the exception that the ambient dimension of $X$ has been increased from $3$ to $6$. This segment of the simulation is designed to explore the persistence of the notable efficacy of our proposed semi-supervised approaches when a significant disparity exists between the intrinsic dimension and the ambient dimension of the feature space. Given that the increase in the dimension of $X$ intensifies the challenge of model training, we opt to raise the size of labeled data to $n = 200$. And the unlabeled data size $m$ is selected from $\{0, 200, 500, 1000, 2000, 3000\}$. 
	
\begin{figure}[ht!]
		\centering
		\subfloat[setting I under snr=2]{\includegraphics[width=0.5\linewidth]{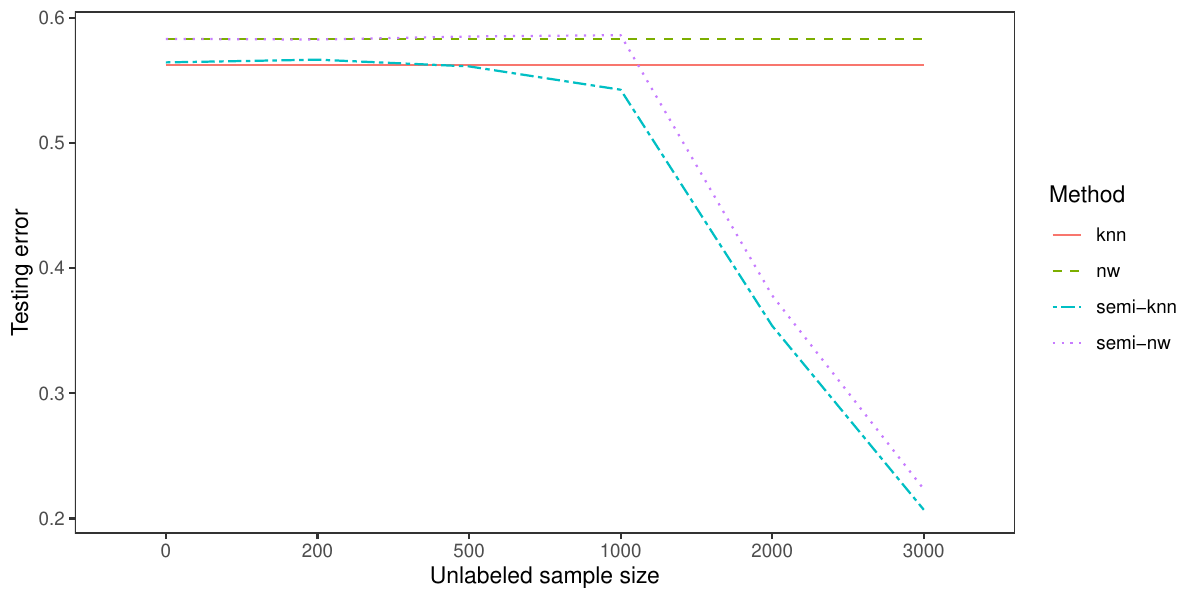}}
		\subfloat[setting I under snr=4]{\includegraphics[width=0.5\linewidth]{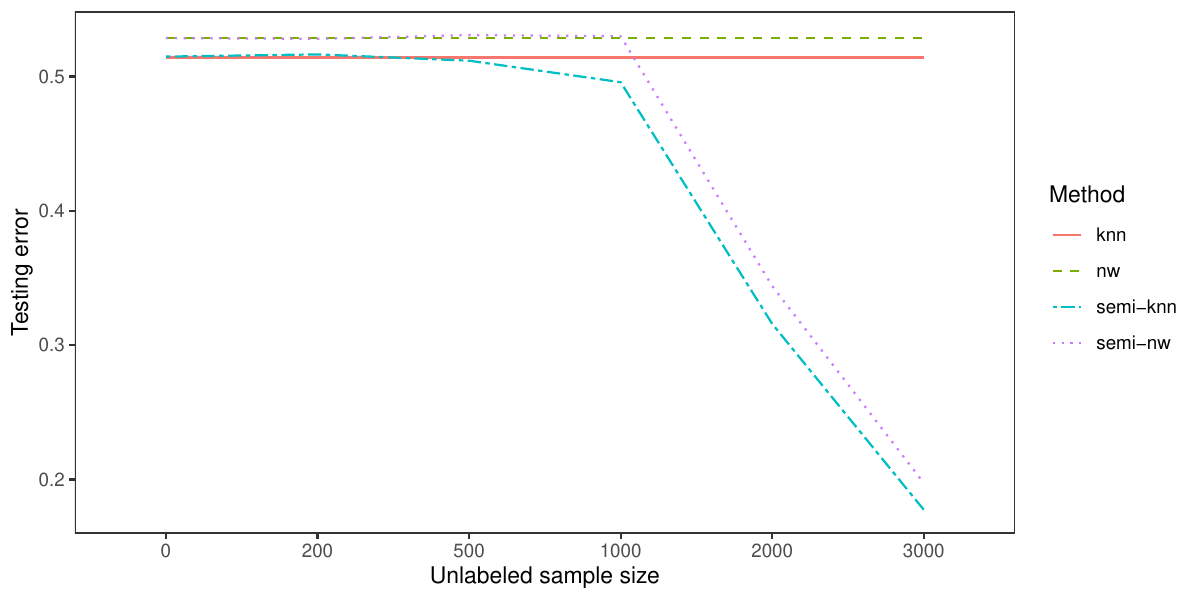}}
		\vspace{4mm}
		\subfloat[setting II under snr=2]{\includegraphics[width=0.5\linewidth]{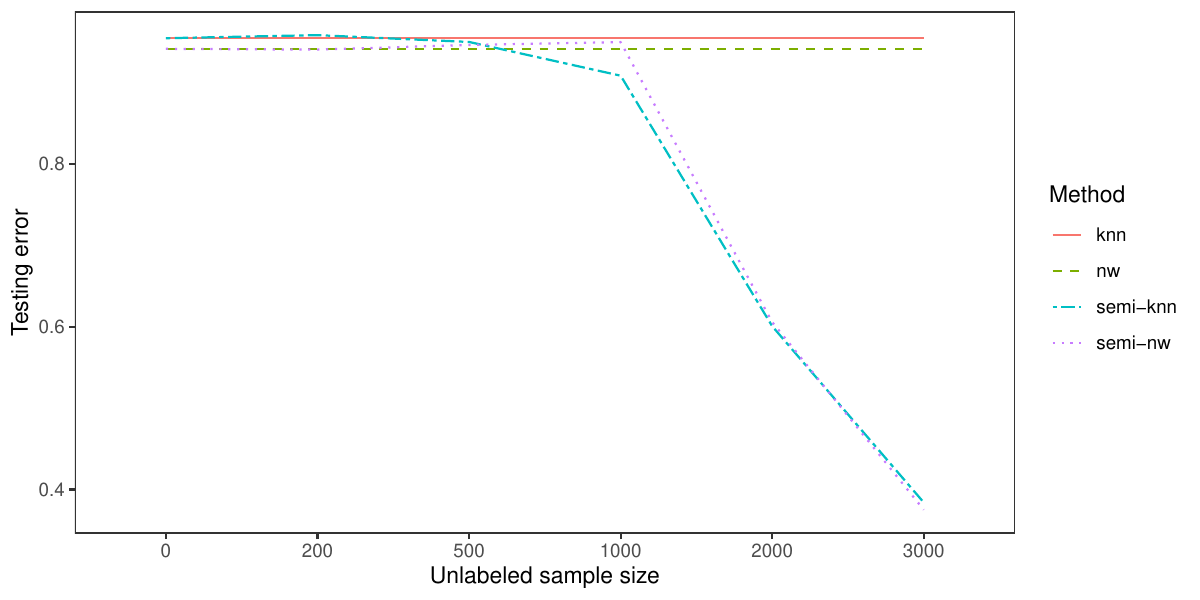}}
		\subfloat[setting II under snr=4]{\includegraphics[width=0.5\linewidth]{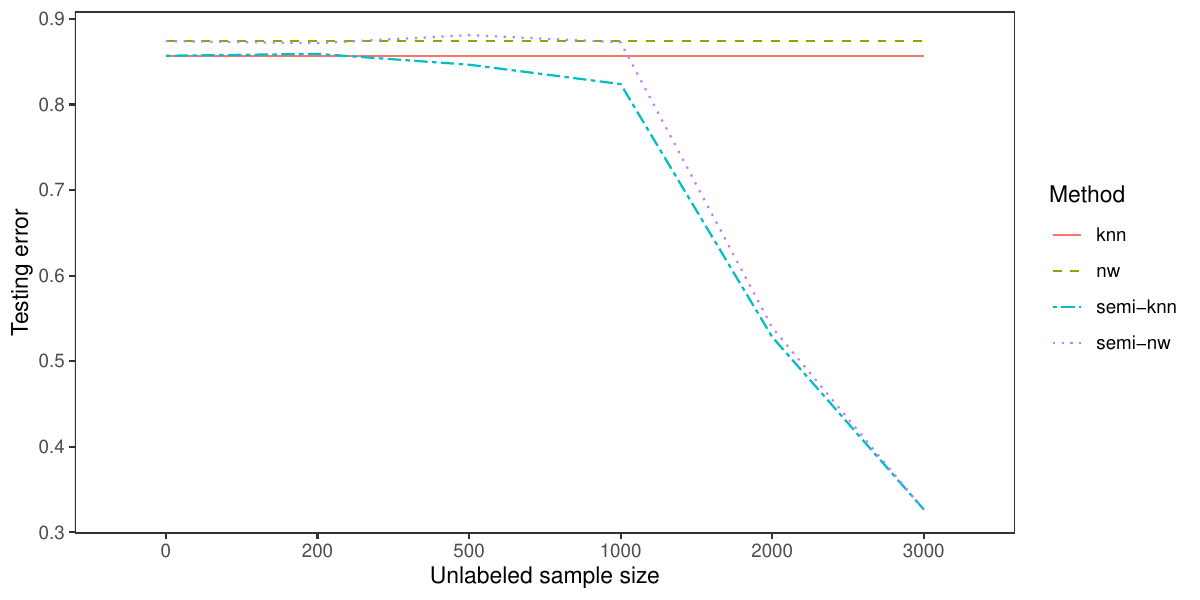}}
		\caption{AMSE of different methods for setting I, II with $200$ labeled points under snr$=2, 4$ when the dimension of $X$ is equal to $6$.}
		\label{fig:spd_six dimension}
	\end{figure}

	\emph{4.2.1. Responses for symmetric positive-definite matrices.}
	Similarly, we generate $Y$ using setting I, II as described in Section 4.1.1 and consider two distinct noise levels, snr$=2$ or $4$. However, unlike the previous simulations, we now delve into the Fr\'echet regression relationship between $Y$ and $6$-dimensional $X$. The results presented in Figure~\ref{fig:spd_six dimension} reveal that a large amount of unlabeled data still brings notable benefits to our semi-supervised predictions, reaffirming the efficacy of our methods. Nonetheless, in contrast to the situation with $3$-dimensional $X$, it is apparent that more unlabeled samples is required to distinctly differentiate the performance of the semi-supervised and supervised methods. Additionally, with an increased number of labeled samples (in contrast to the previous labeled sample size setting of $n = 100$), the kNN methods no longer perform worse than the NW methods. Remarkably, with more labeled samples, the predictive performance of supervised methods is slightly degraded compared to the previous simulations. It could be explained by the curse of dimensionality. However, the prediction accuracy of the semi-supervised methods can be further improved.  For example, when the size of unlabeled data is $m=3000$, a comparative observation of (c) or (d) within Figure~\ref{fig: spd_three dimension} and \ref{fig:spd_six dimension} tells that the AMSE is reduced for both semi-supervised kNN Fr\'echet regression and semi-supervised NW Fr\'echet regression. This fact somewhat implies to us that the performance of our semi-supervised methods depends on the intrinsic dimension of $X$ when there are enough unlabeled data, which is also in line with the theoretical guarantee.

\begin{figure}[ht!]
		\centering
		\subfloat[setting III]{\includegraphics[width=0.5\linewidth]{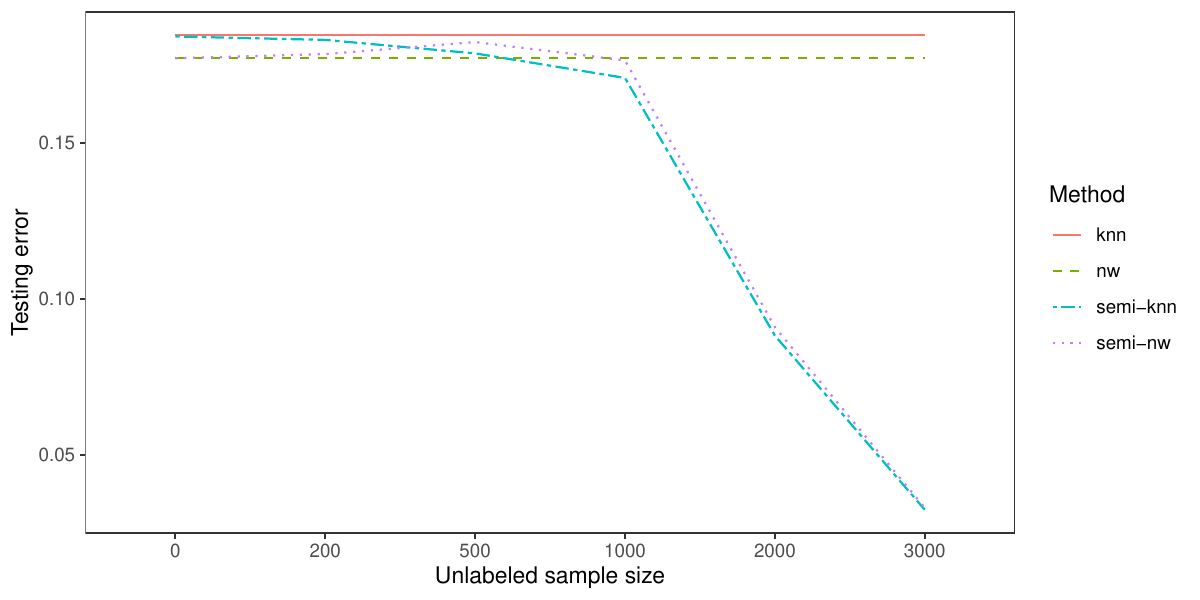}}
		\subfloat[setting IV]{\includegraphics[width=0.5\linewidth]{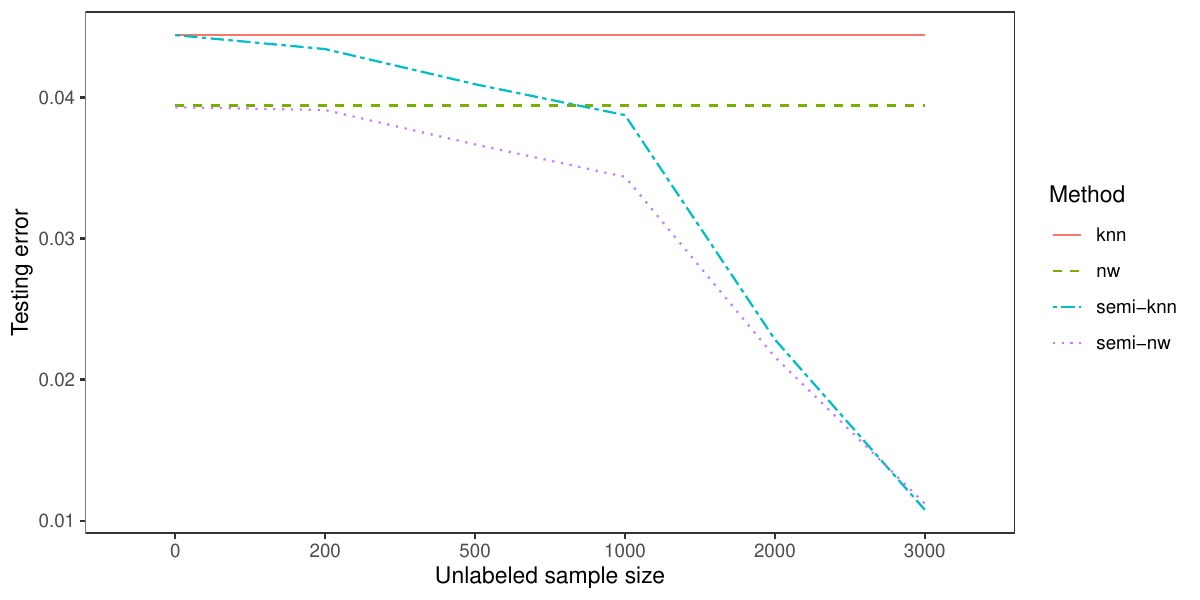}}
		\caption{AMSE of different methods for setting III, IV with $200$ labeled points when the dimension of $X$ is equal to $6$.}
		\label{fig:sph_six dimension}
	\end{figure}

	\emph{4.2.2. Responses for spherical data.}
	Then we turn to regress the spherical data on $6$-dimensional $X$ under the settings III and IV given in Section 4.1.2.  The results are shown in Figure~\ref{fig:sph_six dimension}. The analysis here is akin to Section 4.2.1 and will not be reiterated.
	
\begin{figure}[ht!]
		\centering
		\subfloat[setting I under snr=4]{\includegraphics[width=0.5\linewidth]{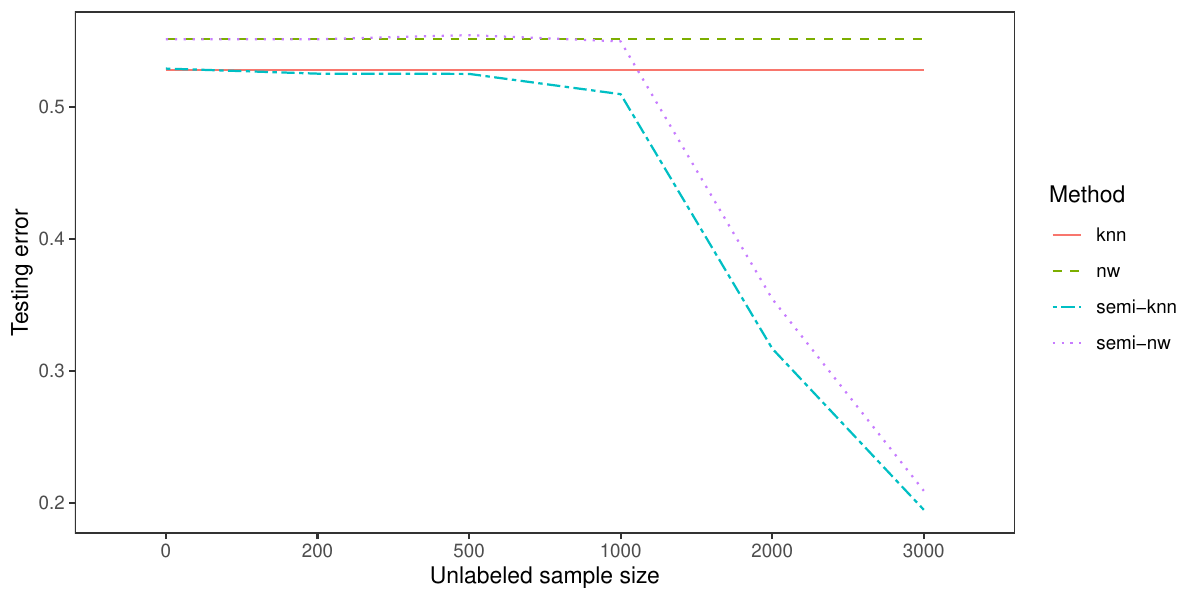}}
		\subfloat[setting III]{\includegraphics[width=0.5\linewidth]{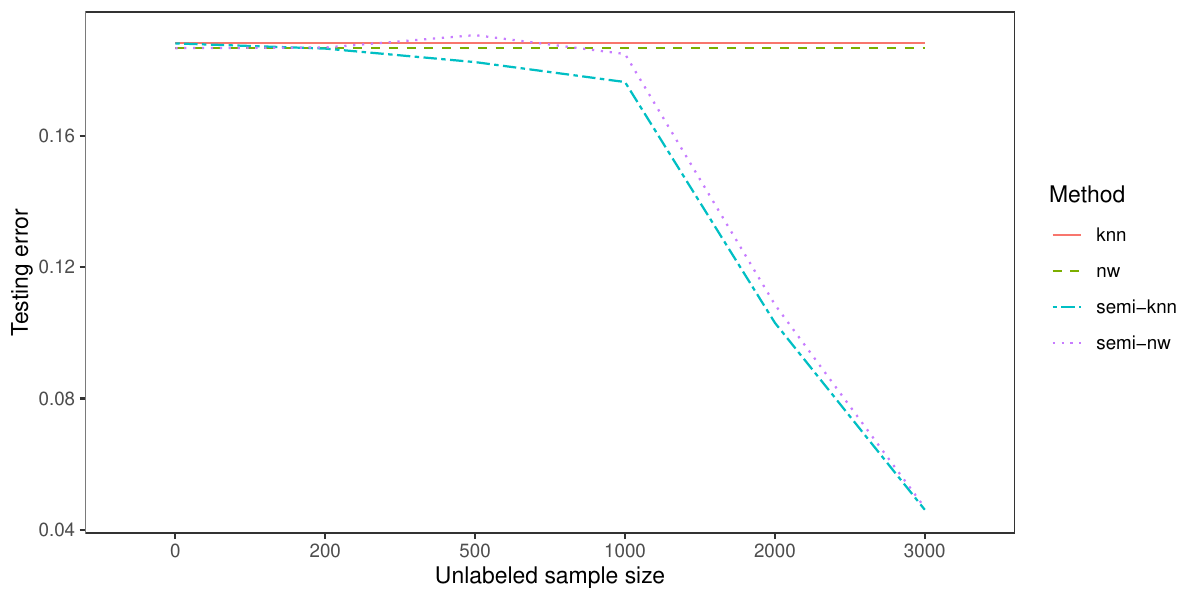}}
		\caption{AMSE of different methods for setting I, III with $200$ labeled points when the dimension of $X$ is equal to $6$ and the components of $U$ are correlated.}
		\label{fig:sph_six dimension_cor}
	\end{figure}

	Lastly, we revisit the case where $U$ is drawn from the aforementioned truncated multivariate normal distribution, and the results are depicted in Figure~\ref{fig:sph_six dimension_cor}. Again, the behavior of all the methods does not differ much from the case where $U$ comes from a uniform distribution.
	
	\section{Real data}\label{Real data}
\begin{figure}[ht!]
		\centering
		\includegraphics[width=0.6\linewidth]{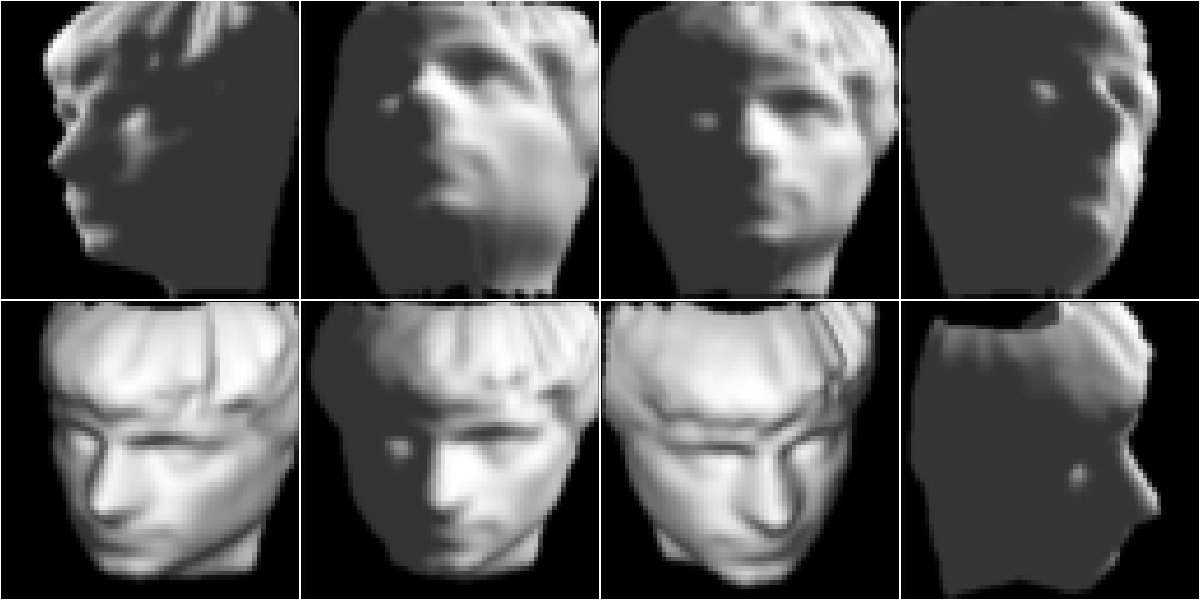} \hspace{0.5cm}
		\includegraphics[width=0.3\linewidth]{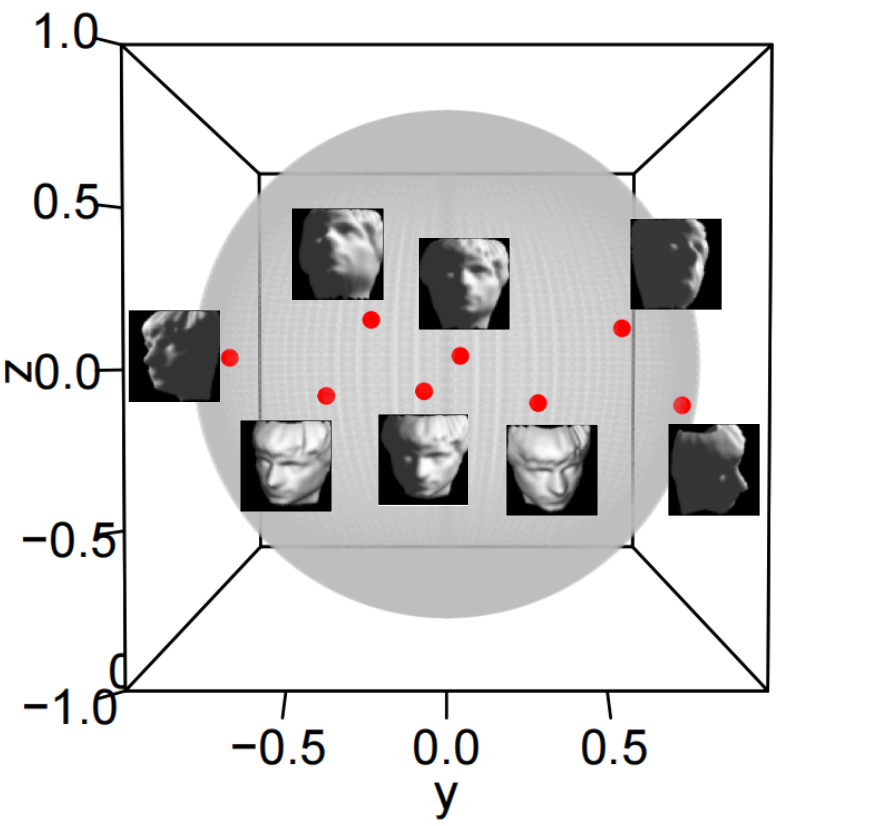}
		\caption{Eight representative samples of face pose images (left) and corresponding points on the sphere (right).}
		\label{fig:face}
	\end{figure}

	In this section, we conduct an analysis on a real dataset of faces. The data are collected from  \url{https://web.archive.org/web/20040411051530/http://isomap.stanford.edu/}.  It comprises $698$ grayscale images of size $64 \times 64$ capturing the face of a single individual, each taken from different angles and light directions. Nevertheless, the intrinsic dimension of the image is 3, as it is determined by three specific directions: the left-right pose angle $(-75^\circ \sim 75^\circ)$, the up-down pose angle $(-10^\circ \sim 10^\circ)$, and the lighting direction. Here we amalgamate the left-right angle and the up-down angle into a singular direction, effectively encapsulating the orientation of the face. This composited direction can be mathematically depicted as a point on the unit sphere, as illustrated in Figure~\ref{fig:face}. We utilize this spherical representation as the target for our Fr\'echet regression, with the corresponding image as the input feature. We proceed by randomly selecting  $70$ or $140$ samples from the dataset as labeled samples, leaving the remaining as unlabeled samples. Subsequently, we apply kNN Fr\'echet regression, semi-supervised kNN Fr\'echet regression, NW Fr\'echet regression and semi-supervised NW Fr\'echet regression once more to predict the response values of the unlabeled samples. For two semi-supervised methods, we construct 4nn-graphs to facilitate the computation of the shortest graph distances. The selection of the kernel function and other hyperparameters follows the choices outlined in Section~\ref{Simulation}. The entire procedure is run $100$ times, and the AMSE with respect to the spherical geodetic distance is used as an evaluation criterion for the merit of all methods.
	
	\begin{figure}[ht!]
		\centering
		\includegraphics[width=1\linewidth]{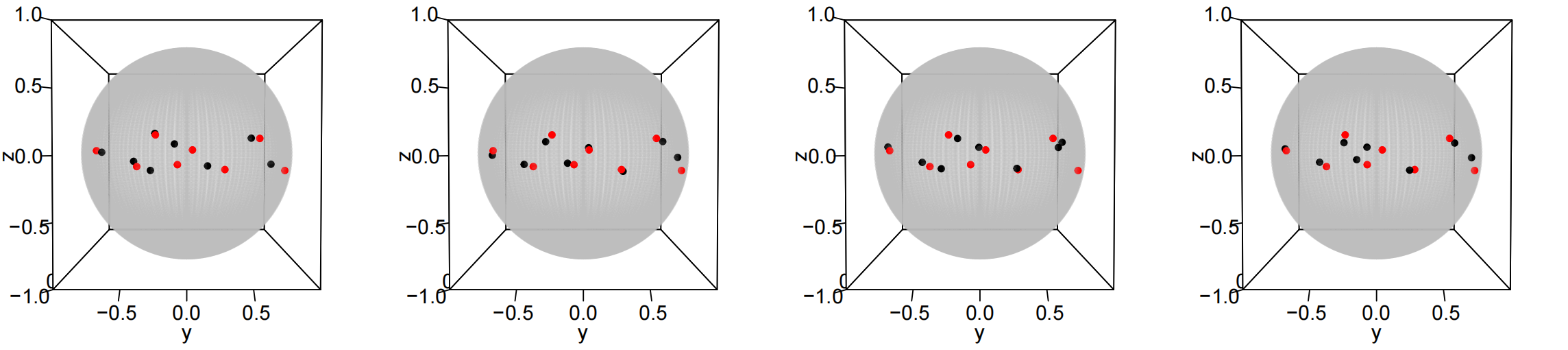} 
		\caption{Predictions of knn, semi-knn, nw and semi-nw (from left to right) for the eight representative images with $70$ labeled data. The red points represent the real directions, and the black points represent the predicted directions.}
		\label{fig:face_methods}
	\end{figure}
	
\begin{figure}[ht!]
		\centering
		\subfloat{
			\begin{minipage}[b]{1\textwidth}
				\centering
				\includegraphics[width=0.8\linewidth]{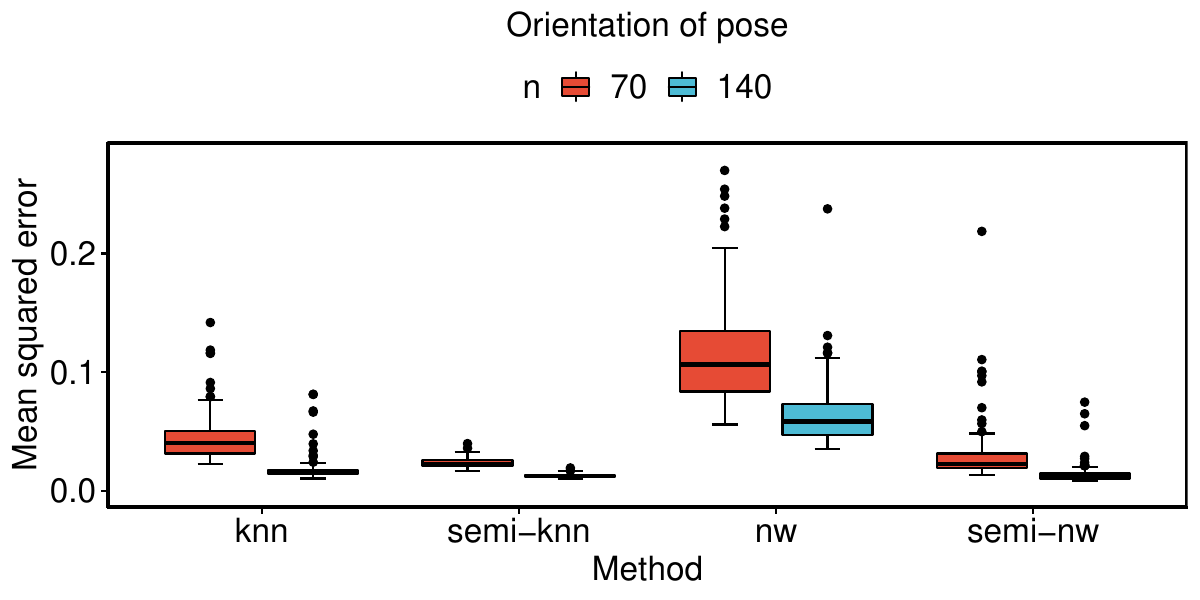}
			\end{minipage}
		}\\
		\subfloat{
			\begin{minipage}[b]{1\textwidth}
				\centering
				\includegraphics[width=0.8\linewidth]{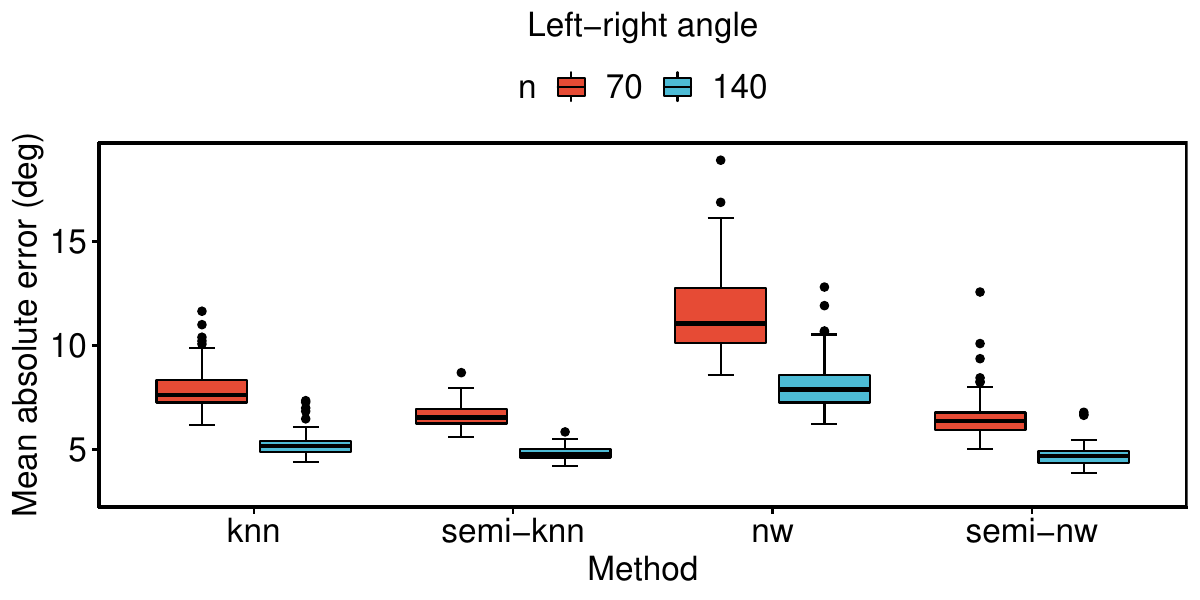}
			\end{minipage}
		}\\
		\subfloat{
			\begin{minipage}[b]{1\textwidth}
				\centering
				\includegraphics[width=0.8\linewidth]{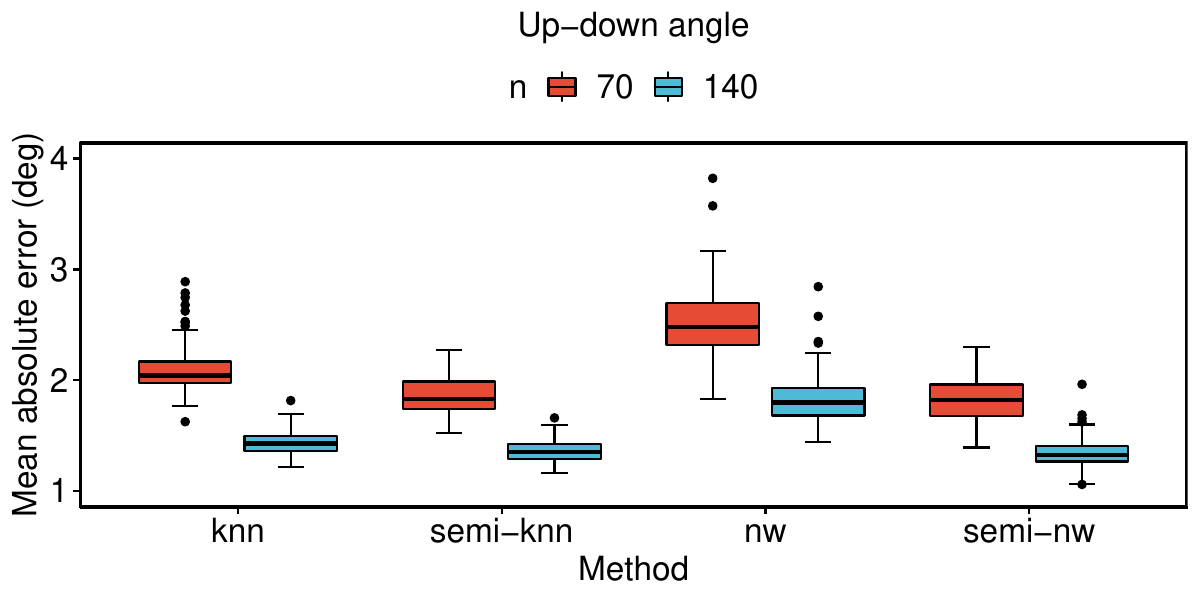}
			\end{minipage}
		}
		\caption{Mean squared errors for geodesic distance, mean absolute error for horizontal and vertical angles using different methods over $100$ realizations with $70$ or $140$ labeled points.}
		\label{real data}
	\end{figure}

	When the size of labeled data is $n=70$, the AMSE is $0.046$ for kNN Fr\'echet regression, $0.024$ for semi-supervised kNN Fr\'echet regression, $0.118$ for NW Fr\'echet regression and $0.031$ for semi-supervised NW Fr\'echet regression. Semi-supervised kNN Fr\'echet regression has the best performance. Clearly, we can derive the corresponding predicted values for the horizontal and vertical orientations of the face pose based on the predicted response values. Then we calculate the average mean absolute error with respect to the angle predictions over the $100$ realizations. The horizontal angle errors for the four methods are $7.886, 6.631, 11.590, 6.567$ respectively. And the vertical angle errors for them are $2.105, 1.857, 2.527, 1.815$ respectively. 
	Instead, the semi-supervised NW Fr\'echet regression makes the most accurate predictions in the single orientation. Hence, the comparative assessment of the merits among different methods is not fixed; it also depends on various evaluation criteria. When the size of labeled data is increased to $n=140$, the three kinds of error can be computed again. The detailed results are recorded as box plots presented in Figure~\ref{real data}. It can be clearly seen that the prediction accuracy of all methods is improved with more labeled data.
	
	To vividly portray the predictions made by each method regarding facial orientation, we select eight highly distinctive images as representatives (refer to Figure~\ref{fig:face}) and plot the predictions of different methods on four unit spheres (see Figure~\ref{fig:face_methods}). Upon visual inspection, it becomes evident that the predictions of semi-supervised kNN Fr\'echet regression closely align with the true values, followed by the semi-supervised NW Fr\'echet regression. In contrast, there are considerable deviations in the predictions of the kNN Fr\'echet regression and NW Fr\'echet regression for certain images. Overall, the semi-supervised methods continue to show significant advantages in analyzing this real data.

	\section{Discussion}\label{Discussion}
	In the realm of semi-supervised regression with Euclidean responses, recent studies have progressively gravitated towards deep learning frameworks. Neural networks have demonstrated significant potential in traditional semi-supervised learning approaches due to their exceptional capacity for feature extraction. However, the existing related literature often lacks comprehensive theoretical guarantees and model interpretations. As for Fr\'echet regression, current research primarily focuses on classical statistical and machine learning models, with no existing development of deep learning variations. This paper addresses the Fr\'echet regression problem in a semi-supervised scenario for the first time. The motivation behind this research stems from two main considerations. Firstly, it acknowledges the high cost associated with obtaining non-Euclidean labels, highlighting the importance of exploring alternative approaches such as semi-supervised learning. Secondly, the article aims to serve as an initial stepping stone for future research in this direction, envisioning the further development and refinement of semi-supervised Fr\'echet regression as research on Fr\'echet regression progresses.
	
	An interesting extension of the current paper is the online semi-supervised regression. Within the inductive learning framework articulated here, new sample points are added by turns to the graph initially constructed upon the entire training samples for prediction purposes. As one might intuitively surmise, with the accumulation of subsequent feature samples, they can serve to augment the pool of unlabeled instances, thus further improving the accuracy of graph distance approximation for the geodesic distance inherent in a low-dimensional manifold. In other words, the continuous introduction of new samples from the feature space can engender a concomitant expansion of the hitherto constructed $r$-graph. Of course, as the sample size increases, the connectivity radius $r$ of the $r$-graph warrants a judicious update. We leave it for future research.

	\bigskip
	
	\bibliographystyle{apalike}
	
	\bibliography{bibfile}
\end{document}